\documentclass[a4paper]{amsart}
\usepackage{graphicx}

\numberwithin{equation}{section}
\newtheorem{thm}[equation]{Theorem}
\newtheorem{pro}[equation]{Proposition}
\newtheorem{cor}[equation]{Corollary}
\newtheorem{lem}[equation]{Lemma}

\theoremstyle{definition}
\newtheorem{defi}[equation]{Definition}
\newtheorem{exa}[equation]{Example}

\theoremstyle{remark}

\newcommand{\End}{\operatorname{End}}

\newcommand{\zz}{{\mathbb{Z}}}
\newcommand{\cc}{{\mathbb{C}}}

\newcommand{\id}{{\operatorname{id}}}

\newcommand{\ve}{{\varepsilon}}
\newcommand{\hec}{{\mathcal{H}_k(q^2)}}

\newcommand{\ur}{{\mathcal{R}}}
\newcommand{\hr}{{\check{R}}}
\newcommand{\uq}{{\mathfrak{U}_q}}
\newcommand{\uu}{{\mathfrak{U}_q(gl(m,n))}}
\newcommand{\mh}{{\mathfrak{h}}}

\newcommand{\mF}{{\mathbb{F}}}
\newcommand{\mg}{{\mathfrak{g}}}

\allowdisplaybreaks

\begin{document}

\title
[Highest weight vectors of irreducible representations of
$\uq\left(gl(m,n) \right)$] {Highest weight vectors of irreducible
representations of the quantum superalgebra $\uq\left(gl(m,n)
\right)$}

\author{
Dongho Moon
}
\address{
School of Mathematics \\
Korea Institute for Advanced Study\\
Seoul, Korea \\
\textrm{ and }
Department of Mathematics \\
Sejong University \\
Seoul, Korea}
\email{dhmoon@sejong.ac.kr}
\thanks{Research supported in part by National Science Foundation grant
DMS-9622447}


\begin{abstract}
The Iwahori-Hecke algebra $\hec$ of type A acts on tensor product
space $V^{\otimes k}$ of the natural representation of the quantum
superalgebra $\uq(gl(m,n))$. We show this action of $\hec$ and the
action of $\uq(gl(m,n))$ on the same space determine commuting
actions of each other. Together with this result and Gyoja's
$q$-analogue of the Young symmetrizer, we construct a highest
weight vector of each irreducible summmand of the tensor product
space $V^{\otimes k}$, for $k=1,2,\ldots $.
\end{abstract}
\maketitle

\setcounter{section}{-1}

\section{Introduction}

One of the main studies of the representation theory of a
semisimple Lie algebra $\mg$ is constructing all the irreducible
$\mg$-modules. Related with this problem, we are also interested
in obtaining highest weight vectors of the irreducible summands of
a tensor representation.

When $\mg$ is the \emph{special linear Lie algebra} $sl(n)$ or the
\emph{general linear Lie algebra} $gl(n)$ over the field $\cc$,
this problem was successfully solved by I.~Schur in \cite{schur1}
and \cite{schur2}. Schur investigated the tensor product space of
the \emph{natural representation}, which is the irreducible
representation of $gl(n)$ with highest weight $\epsilon_1$. He
showed the action of $gl(n)$ on the tensor product space generate
the full centralizer of the symmetric group action. And then, from
the \emph{double centralizer theorem}, we may show the associative
algebra generated by actions commuting with actions of $gl(n)$,
which is called the \emph{centralizer algebra} of $gl(n)$, is a
quotient of the group algebra $\cc S_k$ of the symmetric group
$S_k$. This result is often called \emph{Schur-Weyl duality}, and
it is important for understanding the representation theory of
$gl(n)$. Schur used results  on the representation theory of the
symmetric group $S_k$ by F.~Frobenius \cite{frob} and by A.~Young
\cite{young}. Schur used the decomposition of the group algebra
$\cc S_k$ to obtain the irreducible decomposition of the tensor
product space via the \emph{Young symmetrizers}.

Same approach was made by A.~Berele and A.~Regev \cite{br}  and
G.Benkart and C.~Lee~Shader \cite{nova} for the \emph{general
linear Lie superalgebra} $gl(m,n)$. When $\mg = gl(m,n)$, the
centralizer algebra is again a homomorphic image of $\cc S_k$, and
we can also use the Young symmetrizers to decompose the tensor
product space.

In 1986, M.~Jimbo \cite{jimbo} constructed the
\emph{Drinfel'd-Jimbo} quantum group $\uq(gl(n))$ of $gl(n)$. He
also showed the action of the Iwahori-Hecke algebra of Type A,
$\hec$, on the $k$-fold tensor product space of the natural
representation commutes with the action of $\uq(gl(n))$. And a
$q$-analogue of the Young symmetrizers was obtained by A.~Gyoja
\cite{gyoja}.

The usual trick for proving that the action of general linear Lie
algebra $gl(n)$ generates the full centralizer of the symmetric
group action uses the idempotent $\sum\limits_{\sigma \in S_k}
\sigma$ to construct a projection map onto the $gl(n)$-invariants.
Unfortunately this method is no longer useful for the quantum
case. To show the action of $\hec$ determines the full centralizer
of $\uq(gl(n))$, R.~Leduc and A.~Ram used the path algebra
approach to the centralizer algebra of $\uq(gl(n))$ in  \cite{lr}.
But their approach requires that the tensor product space is a
completely reducible $\mg$-modules.

Recently  in \cite{bkk}, G.~Benkart, S.-J.~Kang, and M.~Kashiwara
showed the completely reducibility of the tensor product space of
the natural representation of $\uq(gl(m,n))$ using crystal bases.
Now we may use their result to obtain the full centralizer of
$\uq(gl(m,n))$ in Section~\ref{sec:20}.

Our main result will appear in Section~\ref{sec:maximal}. We will
decompose the tensor product space to obtain finite dimensional
irreducible representations of $\uq(gl(m,n))$ using Gyoja's
$q$-analogue of the Young symmetrizers. Finally we construct a
highest weight vector of each irreducible representation using the
$q$-analogue of the Young symmetrizers.

\section{The quantum superalgebra $\uq(gl(m,n))$}\label{sec:10}

We begin by giving a description of the \emph{general linear Lie
superalgebra} $gl(m,n)$.

Let $E_{i,j} \in gl(m,n)$ denote the matrix unit which has
$1$ at $(i,j)$-position and $0$'s at other positions. The Cartan
subalgebra $\mh$ of $gl(m,n)$ is the set of all diagonal matrices
in $gl(m,n)$, which is the $\cc$-span of $H_i$, for $i=1,\ldots,
m+n-1$,  and $J$, where
\begin{equation} \label{eq:cartan}
\begin{aligned}
H_i&= \begin{cases} E_{i, i} - E_{i+1, i+1} &\text{if $i\ne m$,}  \\
                   E_{m,m}+ E_{m+1,m+1} &\text{if $i=m$.}
      \end{cases} \\
  J&= \sum_{k=1}^{m+n} E_{i, i}.
\end{aligned}
\end{equation}

Relative to the adjoint action of the Cartan subalgebra $\mh$,
$gl(m,n)$ decomposes into root spaces
$$
gl(m,n) = \mh \bigoplus \sum_{\alpha \in \mh^*} gl(m,n)_\alpha.
$$

Let
$\epsilon_1,\ldots,\epsilon_m$ and $\delta_1, \ldots,\delta_n$ be
orthonormal bases of $\mathbb{R}^m$ and $\mathbb{R}^n$
respectively. We also write $\epsilon_{m+1}:=\delta_1, \ldots,
\epsilon_{m+n}:=\delta_n$.   The simple roots $\alpha_i \in \mh^*$, and
the fundamental weights $\omega_i\in \mh^*$, for $i=1,\ldots, m+n-1$, are
given by
\begin{align}
\alpha_i&= \begin{cases}
   \epsilon_i-\epsilon_{i+1} \quad  &1\le i \le m-1, \\
   \epsilon_m-\delta_1 \quad &\text{for $i=m$ }, \\
   \delta_{i-m}-\delta_{i+1-m}, \quad &m+1\le i \le m+n-1,
   \end{cases} \label{eq:root} \\
\omega_i&=\epsilon_1+\cdots \epsilon_i \label{eq:weight}.
\end{align}

Let $P$ be the $\zz$-span of $\{\epsilon_1,\ldots,
\epsilon_{m+n}\}$, which we call the \emph{lattice of integral
weights}. And the \emph{dual weight lattice} $P^\vee\subset \mh$
is the free $\zz$-lattice spanned by $E_{i,i}$, $1\le i\le m+n$.

We may define the value $\lambda(h)$ for any $\lambda \in \mh^*$
and $h\in \mh$ by
$$
\epsilon_i(E_{j,j})= \delta_{ij}.
$$
and extending it by linearity. This allows us to define a natural
pairing  $\langle \, , \, \rangle$ between $\mh$ and $\mh^*$ so
that
$$
\langle H_i, \alpha_j \rangle= \alpha_j(H_i)=a_{ij}.
$$

The \emph{Cartan} matrix $A=\left( a_{ij} \right)_{1\le i,j \le
m+n-1}$, where $a_{ij} = \alpha_j(H_i)$, of $sl(m,n)$ {\it (or
$gl(m,n)$)} satisfies
\begin{equation} \label{eq:cargl}
a_{ij}=\begin{cases}
2 & \text{ if $i=j \ne m$}, \\
0 & \text{ if $i=j=m$},\\
-1 & \text{ $j= i-1$ or  $j=i+1, i\ne m$}, \\
1 & \text{ $i=m$, $j=m+1$}, \\
0 & \text{ otherwise}.
\end{cases}
\end{equation}
Note that the Cartan matrix $A$ is symmetrizable, i.e., if we
define $d_i$  for  $i\in I$ by
\begin{equation} \label{eq:di}
d_i= \begin{cases}
1, &\text{ if $1 \le i \le m$}, \\
-1, &\text{ if $m+1\le i \le m+n-1$},
\end{cases}
\end{equation}
and if we let $D=\text{diag}(d_1,\ldots,d_{m+n-1})$ be the
the diagonal matrix with diagonal entries $d_i$ for $i\in I$,
then $A^{\text sym}=DA$ is a symmetric matrix.

A root $\alpha$ is \emph{even} if $gl(m,n)_\alpha \cap
gl(m,n)_{\bar 0} \ne \{0\}$ or \emph{odd} if $gl(m,n)_\alpha \cap
gl(m,n)_{\bar 1} \ne \{0\}$. Hence all the simple roots except
$\alpha_m$ are even.

Now we give a definition of the quantum superalgebra
$\uq(gl(m,n))$. A classical contragredient Lie superalgebra of
rank $r$ can be characterized by its Cartan matrix
$A=\bigl(a_{ij}\bigr)_{i,j\in I}$ and a subset $\tau \subset
I=\{1,2,\ldots,r\}$ for the \emph{odd} simple roots. For the
general linear Lie superalgebra $gl(m,n)$ of rank $r=m+n-1$, we
have $\tau=\{m\}$.

The Serre-type presentation of $sl(m,n)$ (or $gl(m,n)$) and the
definition of the quantum superalgebra of $sl(m,n)$ (or $gl(m,n)$)
were obtained by various authors all roughly about the same time
(see for example \cite{kt}, \cite{flv}, \cite{scheu92} or
\cite{scheu93}). Readers may refer those papers for the
presentation of $sl(m,n)$.

But in this paper we use a slightly different definition of the
quantum super algebra $\uq(gl(m,n))$ which was appeared in
\cite{bkk} to quote results in the paper. Let $q$ be an
indeterminate and let $\cc(q)$ denote the field of rational
functions in $q$. Let $q_i:=q^{d_i}$, where $d_i$ is defined in
\eqref{eq:di}.

\begin{defi}[\cite{bkk}]
The associated quantized enveloping algebra $\uq(gl(m,n))$ is the
unital associative algebra over $\cc(q)$ with generators $E_i$,
$F_i$ ($i\in I$), $q^h$ ($h\in P^\vee$), which satisfy the
following defining relations:
\begin{align}
 &q^h=1 \quad \text{ for $h=0$}, \\
 &q^{h_1+h_2}= q^{h_1} q^{h_2} \quad \text{ for $h_1,h_2\in
 \mathfrak{h}^*$}, \\
 &q^h E_i = q^{\langle h, \alpha_i \rangle} E_i q^h ,   \\
 &q^h F_i = q^{-\langle h, \alpha_i \rangle} F_i q^h , \quad
 \text{ for $h\in \mathfrak{h}^*$ and $i\in I$}, \\
 &E_iF_j -(-1)^{p(E_i)p(F_j)}F_j E_i = \delta_{i,j} \frac{k_i
 -k_i^{-1}}{q_i - q_i^{-1}}  \quad \text{ for $i,j \in I$, where
 $k_i=q^{d_i H_i}$}, \\
 &E_iE_j -(-1)^{p(E_i)p(E_j)}E_j E_i = 0,  \\
 &F_iF_j -(-1)^{p(F_i)p(F_j)}F_j F_i = 0, \quad \text{if $|i-j| >
 2$}, \\
 & E_i^2E_j-(q_i+q_i^{-1}) E_iE_jE_i +E_jE_i^2=0,  \\
 & F_i^2F_j-(q_i+q_i^{-1}) F_iF_jF_i +F_jF_i^2=0 \quad \text{if
 $|i-j|=1$ and $i\ne m$},  \\
 & E_m^2=F_m^2=0,  \\
\begin{split}
 &E_mE_{m-1}E_mE_{m+1} + E_mE_{m+1}E_mE{m-1} + E_{m-1}E_mE_{m+1}E_m \\
 &\qquad\qquad + E_{m+1}E_mE_{m-1}E_m -(q+q^{-1})E_mE_{m-1}E_{m+1}E_m =0,
\end{split} \\
\begin{split}
 &F_mF_{m-1}F_mF_{m+1} + F_mF_{m+1}F_mF{m-1} + F_{m-1}F_mF_{m+1}F_m \\
  &\qquad \qquad + F_{m+1}F_mF_{m-1}F_m -(q+q^{-1})F_mF_{m-1}F_{m+1}F_m =0.
\end{split}
\end{align}
The parities are given as $p(q^h)=0$ for all $h\in P^\vee$,
$p(E_i)=p(F_i)=0$ for $i \ne m$, and $p(E_m)=p(F_m)=1$.
\end{defi}

It is also worth to write $K_i:= q^{H_i}$ for each $i\in I$ as a
symbol, even though it has no relation with the meaning of
\emph{``the $H_i$th power of $q$"}. Sometimes we will also use
notations $E_{\alpha_i}:=E_i$, $F_{\alpha_i}:=F_i$, and
$H_{\alpha_i}:=H_i$.

We let $\uu_{\ge 0}$ be the subalgebra of $\uu$ generated by
$E_i$ and $q^h$, $h\in P^\vee$. Similarly let $\uu_{\le 0}$ be the
subalgebra generated by $F_i$ and $q^h$, $h\in P^\vee$.

The Hopf superalgebra structure of $\uq(gl(m,n))$ is given by
comultiplication $\Delta: \uq(gl(m,n)) \longrightarrow \uq(gl(m,n))
\otimes \uq(gl(m,n))$ such that
\begin{equation} \label{eq:comul}
\begin{split}
\Delta(E_i)&=E_i\otimes k_i^{-1} + k_i \otimes E_i,   \\
\Delta(F_i)&=F_i\otimes k_i^{-1} + k_i \otimes F_i, \\
\Delta(q^h)&=q^h \otimes q^h .
\end{split}
\end{equation}
The antipode $S: \uq(gl(m,n)) \longrightarrow \uq(gl(m,n))$ is given by
\begin{equation}
\begin{split}
S(E_i)&=-q_i^{-a_{ii}}E_i, \\
S(F_i)&=-q_i^{a_{ii}}F_i, \\
S(q^h) &= q^{-h},
\end{split}
\end{equation}
and the counit $\ve: \uq(gl(m,n)) \longrightarrow \cc(q)$ by
\begin{equation}
\begin{split}
\varepsilon(E_i)&=\varepsilon(F_i)=0, \\
\varepsilon(q^h)&=1.
\end{split}
\end{equation}

A $\uq(gl(m,n))$-module $M$ is called a \emph{weight module} if it
admits a weight space decomposition
\begin{equation}\label{eq:wtmod}
 M=\bigoplus\limits_{\lambda \in P} M_\lambda,
\end{equation}
where $M_\lambda = \{u\in M \mid q^h=u q^{\langle h,\lambda
\rangle}u, \; \text{for any $h\in P^\vee$}\}$. A weight module $M$
is a \emph{highest weight module} with \emph{highest weight}
$\lambda \in P$ if there exists a unique nonzero vector $v_\lambda
\in V$ up to constant multiples such that
 \begin{enumerate}
\item $M=\uq(gl(m,n))v_\lambda$,
\item $E_i v_\lambda=0$ for all $i\in I$, and
\item $q^h v_\lambda = q^{\lambda(h)} v_\lambda$ for all $h\in P^\vee$.
 \end{enumerate}

The set of dominant integral weights is defined by
 \begin{equation}
\Gamma=\left\{
\lambda=\sum_{i=1}^{m+n} \lambda_i \epsilon_i \mid \lambda_i \in \zz,
\lambda_1 \ge \lambda_2 \ge \cdots \ge \lambda_m, \lambda_{m+1}\ge
\cdots \ge \lambda_{m+n} \right\}.
\end{equation}
Let $\Gamma^+$ be a subset  of $\Gamma$ defined by
 \begin{equation} \label{eq:polwt}
\Gamma^+=\left\{\lambda \in \Gamma \mid
\lambda_i \ge 0, \quad 1\le i \le m+n \right\}.
 \end{equation}

Now we define the \emph{fundamental representation} of $\uu$. Let
$V=V_{\bar 0}\oplus V_{\bar 1}=\cc(q)^{m}\oplus \cc(q)^n$ be a
$\zz_2$-graded vector space of dimension $(m+n)$ over $\cc(q)$.
Let $T=\{t_1,\ldots,t_m\}$ be a basis of $V_{\bar 0}$ and
$U=\{u_1,\ldots,u_n\}$ be a basis of $V_{\bar 1}$ so that the
parities of the basis vectors are given by $p(t_i)=0$ and
$p(u_i)=1$. Sometimes it is convenient to write $b_1:=t_1, \ldots,
b_m:=t_m$, and $b_{m+1}:=u_1, \ldots, b_{m+n} :=u_n$.

The fundamental (super) representation
$(\rho,V)$, $\rho : \uu \longrightarrow \End(V)$,
of $\uu$ is given by setting
\allowdisplaybreaks
\begin{equation} \label{eq:fundrep}
\begin{split}
\rho(E_i)&= E_{i,i+1}, \\
\rho(F_i)&= E_{i+1,i}, \\
\rho(q^h)&=\sum_{i=1}^{m+n}q^{\ve_i(h)} E_{i,i}
\end{split}
\end{equation}
It is easy to see that this representation is in fact the
irreducible highest weight module $V(\epsilon_1)$ with highest
weight $\epsilon_1$.

Also we let $W=W_{\bar0} \oplus W_{\bar 1}= \cc^m \oplus \cc^n$ be
a $\zz_2$-graded vector space over $\cc$. Note that we may suppose
$V_{\bar 0}=\cc^m \otimes_{\cc} \cc(q)$, $V_{\bar 1}=\cc^n
\otimes_{\cc} \cc(q)$, and $V^{\otimes k}=W^{\otimes k}
\otimes_{\cc} \cc(q)$, which has $\zz_2$-grading $(V^{\otimes
k})_{\bar 0} = (W^{\otimes k})_{\bar 0} \otimes_{\cc} \cc(q)$ and
$(V^{\otimes k})_{\bar 1}= (W^{\otimes k})_{\bar 1} \otimes_{\cc}
\cc(q)$. Also we regard $\End(V^{\otimes k}) = \End(W^{\otimes k})
\otimes_\cc \cc(q)$, which is a $\zz_2$-graded representation
similarly. Since $\uu$ is a Hopf superalgebra, the tensor product
representation $(\rho^{\otimes k}, V^{\otimes k})$ of $\rho$ is a
well-defined super representation for each $k\ge 1$. There is also
a representation $\Psi: \cc S_k \longrightarrow \End (W^{\otimes
k})$ of the group algebra $\cc S_k$ given by $\zz_2$-graded place
permutation on simple tensors.

\section{The universal $R$-matrix, the Hecke algebra, and Gyoja's
$q$-analogue of the Young Symmetrizer}\label{sec:20}

In this section we recall the definition of the universal
$R$-matrix of $\uu$ which appeared in \cite{kt} and show that
there is an action of a certain Hecke algebra $\hec$ on
$V^{\otimes k}$ coming from the universal $R$-matrix, which
commutes with the action of $\uu$ on $V^{\otimes k}$.

Let
$$
\theta: \uu\otimes \uu \longrightarrow  \uu\otimes \uu
$$
be given by $\theta(x\otimes y)=(-1)^{p(x)p(y)} y\otimes x$. We
define opposite comultiplication $\Delta'$ by $\Delta'=
\theta\Delta$.
\begin{thm}[Khoroshkin and Tolstoy \cite{kt}]
There is a unique invertible solution $\ur=\sum\limits_i x_i\otimes y_i
\in \widehat{\uu\otimes \uu}$ (the completion of $\uu\otimes \uu$)
of parity $0$ of the equations
\begin{align}
&\Delta'(x)=\ur\Delta(x) \ur^{-1} \qquad \text{ for all $x\in \uu$}, \\
\begin{split}
&(\Delta\otimes \id) \ur = \ur^{13}\ur^{23}, \\
&(\id\otimes \Delta) \ur = \ur^{13}\ur^{12},
\end{split}
\end{align}
where $\ur^{12}=\sum\limits_i x_i\otimes y_i \otimes 1$, $\ur^{23}=
\sum\limits_i 1\otimes x_i\otimes y_i$, and $\ur^{13}=\sum\limits_i x_i
\otimes 1 \otimes y_i$.
\end{thm}

The universal $R$-matrix $\ur$ is given explicitly in \cite{kt}.
Let $R\in \End( V\otimes V)$ be the transformation induced by the
action of $\ur$ on $V\otimes V$.  Applying $\ur$ to $V\otimes V$
relative to the basis $\{b_i\otimes b_j \mid  i,j=1,\ldots,
m+n\}$, we may compute the matrix of $R$ in  $\End (V\otimes V)$,
which is given by
\begin{equation}
\begin{split}
R=\sum_{i=1}^m q^2 E_{i,i}\otimes &E_{i,i} + \sum_{i=m+1}^{m+n}
E_{i,i} \otimes E_{i,i}  \\
&+ \sum_{i\ne j} qE_{i,i}\otimes E_{j,j} + \sum_{i<j}
(-1)^{p(b_i)}(q^2-1) E_{j,i}\otimes E_{i,j}.
\end{split}
\end{equation}

Let $\hr=\sigma R$, where $\sigma:V\otimes V\longrightarrow V\otimes
V$ is given by $\sigma(v\otimes w) = (-1)^{p(v)p(w)} w\otimes v$.
Then
\begin{equation}
\begin{split}
\hr=\sum_{i=1}^m q^2 E_{i,i}\otimes &E_{i,i} - \sum_{i=m+1}^{m+n}
E_{i,i} \otimes E_{i,i}  \\
&+ \sum_{i\ne j} (-1)^{p(b_i)}q E_{j,i}\otimes E_{i,j} + \sum_{i<j}
(q^2-1) E_{i,i}\otimes E_{j,j}.
\end{split}
\end{equation}
Note that, for homogeneous elements $X\otimes Y \in \End(V\otimes
V)=\End(V)\otimes \End(V)$ and $v\otimes w \in V\otimes V$, we
have $(X\otimes Y)(v\otimes w)= (-1)^{p(Y)p(v)} Xv\otimes Yw$. And
also product of tensors is given as $(X_1\otimes X_2) (Y_1\otimes
Y_2)= (-1)^{p(X_2)p(Y_1)} X_1 Y_1\otimes X_2 Y_2$ for $X_1\otimes
X_2, Y_1\otimes Y_2 \in \End(V\otimes V)$.

By direct calculation, we see that
\begin{equation}\label{eq:r2}
\hr^2 +(1-q^2)\hr= q^2 I_{V\otimes V}.
\end{equation}

For each $j=1,\ldots, k-1$, let
\begin{equation*}
r_j=\id_V^{\otimes j-1} \otimes \hr \otimes \id_V^{\otimes k-j-1} \in
\End(V^{\otimes k}),
\end{equation*}
where $\hr$ operates on the $j$th and the $(j+1)$st tensor slots.
Then using arguments similar to those in
\cite[Proposition~2.18]{lr}, we have
\begin{pro}\label{pro:braid}
\begin{itemize}
\item[(1)] Each $r_j$ commutes with the actions of $\uu$ on $(V^{\otimes
k})$. In other words each $r_j$ is in $\End_{\uu}(V^{\otimes k})$.
\item[(2)] The braid relations are satisfied:
\begin{alignat*}{2}
r_ir_j&=r_jr_i, \qquad &&\text{ for $|i-j|\ge 2$}, \\
r_i r_{i+1} r_i&= r_{i+1} r_i r_{i+1}, \qquad &&\text{ for $1\le i
\le m-2$}.
\end{alignat*}
\end{itemize}
\end{pro}

Also from \eqref{eq:r2} we know
\begin{equation} \label{eq:ri2}
(r_i + \id) (r_i - q^2\id) =0.
\end{equation}

\begin{defi}
The Iwahori-Hecke algebra of type A, denoted by $\hec$, is the
associative algebra over $\cc(q)$ generated by $1,h_1,\ldots, h_{k-1}$
subject to the relations
\begin{itemize}
\item[(B1)] $h_ih_j=h_jh_i$, if $|i-j|\le 2$,
\item[(B2)] $h_ih_{i+1} h_i= h_{i+1}h_i h_{i+1}$,  for
$1\le i\le k-2$,
\item[(B3)] $(h_i+1)(h_i - q^2) =0$.
\end{itemize}
\end{defi}

Notice that $\hec$ is a $q$-analogue of the group algebra $\cc
S_k$ of the symmetric group $S_k$ in the sense that when $q$ is
specialized to $1$, $\hec$ is isomorphic to $\cc S_k$.  Let
$\sigma=s_{i_1}s_{i_2}\cdots s_{i_l}$ be a reduced expression for
$\sigma\in S_k$, where $s_j$ is the transposition $(j\, j+1)$,
$j=i_1, \ldots, i_j$. Then we let that $h(\sigma)=h_{i_1} \cdots
h_{i_l}$. This does not depend on the reduced expression of
$\sigma$. So $h(\sigma_1\sigma_2)=h(\sigma_1)h(\sigma_2)$ if and
only if $\ell(\sigma_1\sigma_2)=\ell(\sigma_1)+\ell(\sigma_2)$.

{F}rom \eqref{eq:ri2} and Proposition~\ref{pro:braid},  we see the
following:
\begin{pro}
There is a representation
$$
\Psi_q: \hec \longrightarrow
\End_\uu(V^{\otimes k})
$$
of the Iwahori-Hecke algebra $\hec$ given by $h_i \mapsto r_i$.
\end{pro}
Note that $\Psi_q$ is a $q$-deformation of $\Psi : \cc S_k
\longrightarrow \End_{\mathfrak{U}(gl(m,n)}\left( W^{\otimes
k}\right)$.

Write $\lambda \vdash k$ to denote that $\lambda$ is a  partition
of $k$, and let $\ell(\lambda)$ denote the number of nonzero parts
of $\lambda$. Corresponding to $\lambda\vdash k$ is its {\em Young
frame} having $k$ boxes with $\lambda_i$ boxes in the $i$th row
and with the boxes in each row left justified. We let $\lambda^*$
be the {\em conjugate partition} of $\lambda$ whose frame is
obtained by reflecting that of $\lambda$ about the main diagonal.
Then $\lambda_j^*$ is just the number of boxes in the $j$th column
of $\lambda$. Note that a partition $\lambda \vdash k$ may be
idenfitied with the dominant weight
\begin{equation*}
\lambda_1 \epsilon_1+\cdots +\lambda_m\epsilon_m
+\lambda'_1\delta_1+\cdots +\lambda'_n\delta_n\in \Gamma^+,
\end{equation*}
where
\begin{equation*}
\lambda'_j=\mathrm{max} \{ \lambda^*_j-m,0 \}.
\end{equation*}

A partition $\lambda$ is said to be of  $(m,n)$ {\em hook-shape}
if $\lambda_{m+1} \le n$. We let $H(m,n;k)$ denote the set of all
partitions of $k$ which are  of $(m,n)$ hook-shape.

The irreducible representations of $S_k$ over  any field  $\mF$ of
characteristic $0$ are indexed by the partitions $\lambda \vdash
k$. We have
\begin{equation} \label{eq:decom}
\mF S_k= \bigoplus\limits_{\lambda\vdash k} I_\lambda,
\end{equation}
where $I_\lambda$ is a simple ideal of $\mF S_k$ which is isomorphic
to a matrix algebra $M_{d_\lambda}(\mF)$.
Here $d_\lambda$ is the dimension of the irreducible $S_k$-module
labeled by $\lambda$.

The following is well-known:
\begin{thm}[See for example, Lusztig \cite{lusz81}]\label{thm:heciso}
The Hecke algebra $\hec$ and  the group algebra $\cc(q) S_k$ of
the symmetric group $S_k$ over the field $\cc(q)$ are isomorphic
as associative algebras.
\end{thm}

Therefore irreducible representations of $\hec$ are also indexed
by the partitions $\lambda \vdash k$, and we also have that

\begin{equation}
\hec= \bigoplus\limits_{\lambda\vdash k} I^q_\lambda,
\end{equation}
where $I^q_\lambda$ is a simple ideal of $\hec$ which is
isomorphic to a matrix algebra $M_{d_\lambda}(\cc(q))$.

A $q$-analogue of the Young symmetrizers is obtained  by Gyoja in
\cite{gyoja}. A {\em standard tableau} $T$ of shape $\lambda
\vdash k$ is obtained by filling in the frame of $\lambda$ with
elements of $\{1,\ldots,k\}$, so that the entries increase across
the rows from left to right and down the columns. Associated to
$\lambda$ are two standard tableaux $S_+=S^+_\lambda$ and
$S_-=S^-_\lambda$, which we illustrate by the following example:

\begin{exa}
If
\begin{equation*}
 \lambda =\rotatebox[origin=Bc]{270}{\rotatebox[origin=tr]{90}
 {\scalebox{0.4}{\includegraphics{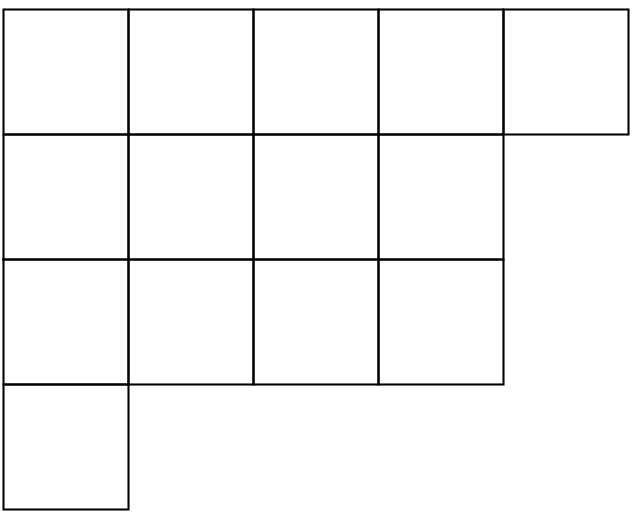}}}},
\end{equation*}
then
\begin{equation*}
S_+=\rotatebox[origin=Bc]{270}{\rotatebox[origin=tr]{90}
{\scalebox{0.4}{\includegraphics{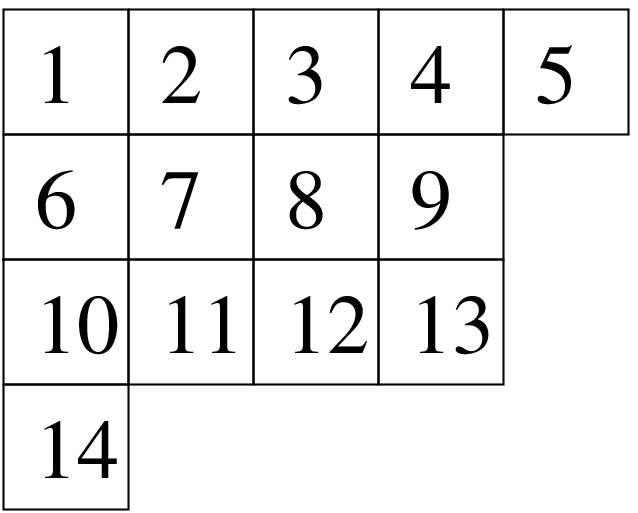}}}}, \quad \text{ and
} S_-=\rotatebox[origin=Bc]{270}{\rotatebox[origin=tr]{90}
{\scalebox{0.4}{\includegraphics{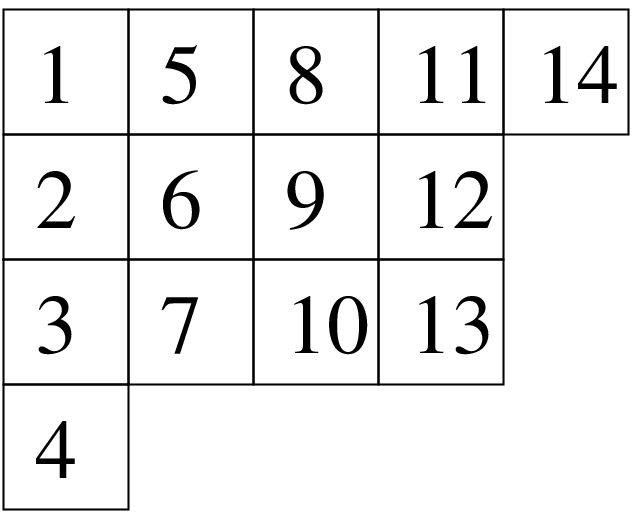}}}}.
\end{equation*}
\end{exa}
Note from the example that the entries of $S_+$ increase by one
across the rows from left to right, and the entries of $S_-$
increase by one down the columns.

Let $T$ be a standard tableau. Let $R(T)$ be the row group of
elements of $S_k$ which permute the entries within each row and
$C(T)$ be the column group of $T$ of permutations which permute
the entries within each column.  Now, for $\lambda\vdash k$, let
\begin{align*}
e_+&=e^+_\lambda:= \sum_{\sigma\in R(S_+)} h(\sigma), \\
e_-&=e^-_\lambda:=\sum_{\sigma\in C(S_-)}(-q^2)^{-\ell(w)}
h(\sigma).
\end{align*}
Then $e_+$ and $e_-$ have the following important properties (see
\cite{gyoja}) :
\begin{alignat}{2}
h(\sigma) e_+&= e_+ h(\sigma) = q^{2 \ell(\sigma)} e_+,
\quad&\text { for $\sigma\in
R(S_+)$}, \label{eq:qmul1}\\
h(\sigma) e_-&= e_- h(\sigma) = (-1)^{\ell(\sigma)} e_-,
\quad&\text{ for $\sigma \in C(S_-)$}.
\end{alignat}

Let $S$ and $T$ be two standard tableaux of shape $\lambda \vdash
k$. We let $\sigma^T_S$ denote the permutation which transforms
$S$ to $T$. We also write $\sigma^T_\pm$ (respectively
$\sigma^\pm_T$, $\sigma^\pm_\mp$) for $\sigma^T_{S_\pm}$
(respectively $\sigma^{S_\pm}_T$, $\sigma^{S_\pm}_{S_\mp}$).  For
example, if
$$
\lambda =\rotatebox[origin=Bc]{270}{\rotatebox[origin=tr]{90}{
\scalebox{0.4}{\includegraphics{5441bas.eps}}}}, \quad \text{ and
} T =\rotatebox[origin=Bc]{270}{\rotatebox[origin=tr]{90}{
\scalebox{0.4}{\includegraphics{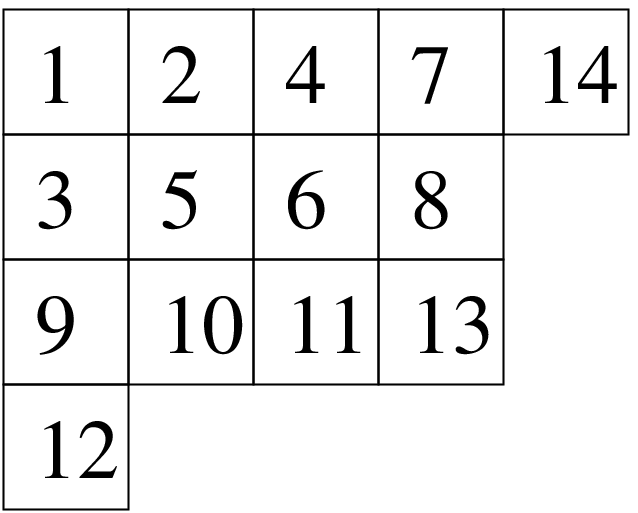}}}},
$$
then
$$
\setcounter{MaxMatrixCols}{14} \sigma^T_+=
\begin{pmatrix}
1&2&3&4& 5&6&7&8&9&10&11&12&13&14 \\
1&2&4&7&14&3&5&6&8& 9&10&11&13&12
\end{pmatrix},
$$
and
$$
\setcounter{MaxMatrixCols}{14} \sigma^T_-=
\begin{pmatrix}
1&2&3&4&5&6&7&8&9&10&11&12&13&14 \\
1&3&9&12&2&5&10&4&6&11&7&8&13&14
\end{pmatrix}.
$$

Let $T$ be a standard tableau of shape $\lambda \vdash k$. Define
$x_T(q) \in \hec$ as
\begin{equation}
x_T(q)= h(\sigma^T_-)
e^-_\lambda\left(h(\sigma^T_-)\right)^{-1}h(\sigma^T_+)
e^+_\lambda \left(h(\sigma^T_+)\right)^{-1}.
\end{equation}
Then there exists a $\xi \in \cc(q)$ depending on the shape
$\lambda$ of $T$ that
$$
x_T(q)x_T(q)= \xi x_T(q).
$$
Now Gyoja's $q$-analogue of the Young symmetrizer is
$$
y_T(q):= \frac{1}{\xi} x_T(q).
$$

Let $T_1$ and $T_2$ be two standard tableaux of same shape. We
compare the entries of $T_1$ and $T_2$ lexicographically starting
at the left end of the first row and moving from left to right. If
all the corresponding entries in the first row are equal, then we
proceed to the second row, etc. So if the first nonzero difference
$j_1 -j_2$ is positive for corresponding entries $j_1$ in $T_1$
and $j_2$ in $T_2$, then we say $T_1 > T_2$.

\begin{pro} [See \cite{gyoja}] \label{pro:gyoja}
The set of all $y_T(q)$'s is a set of primitive idempotents in the
Iwahori-Hecke algebra $\hec$, and so
\begin{itemize}
\item [(1)] For each tableau $T$, $y_T(q)y_T(q)=y_T(q)$.
\item [(2)] If $T_1$  has shape $\lambda\vdash k$ and $T_2$ has shape $\mu
\vdash k$ and $\lambda \ne \mu$,  $y_{T_1}(q)y_{T_2}(q)=0$. And if
$T_1$ and $T_2$ are of the same shape and $T_1 < T_2$, then
$y_{T_1}(q) y_{T_2}(q)=0$.
\end{itemize}
\end{pro}
Note when $q \rightarrow 1$, then $y_T(q)$ specializes to the
Young symmetrizer $y_T$ corresponding to $T$ in the standard case.

\section{Centralizer Theorem}
In this section we show that the actions of $\hec$ determine the
full centralizer of $\uu$.

{F}rom the work by Berele and Regev, we know
\begin{thm}[See \cite{br}] The image $\Psi(\cc S_k)$ is given by
\begin{align*}
\Psi (\cc S_k)\cong & \cc S_k \left/
\Big( \bigoplus_{\substack{\lambda \vdash k \\
\lambda \notin H(m,n)}} I_\lambda \Big )\right.\\
\cong & \bigoplus_{\substack{ \lambda \vdash k \\ \lambda\in H(m,n)}}
I_\lambda.
\end{align*}
\end{thm}

The representation $\Psi_q$ of the Iwahori-Hecke algebra $\hec$ is
completely reducible because $\hec$ is semisimple by
Theorem~\ref{thm:heciso}. Moreover we have
\begin{cor} The image $\Psi_q(\hec)$ is given by
\begin{align*}
\Psi_q (\hec)\cong & \hec \left/
\Big( \bigoplus_{\substack{\lambda \vdash k \\
\lambda \notin H(m,n)}} I^q_\lambda \Big)\right. \\
\cong & \bigoplus_{\substack{ \lambda \vdash k \\ \lambda\in H(m,n)}}
I^q_\lambda.
\end{align*}
\end{cor}
\begin{proof}
Because $\Psi_q$ is a $q$-deformation of $\Psi$, we have
\begin{equation}\label{eq:ee}
\begin{split}
 (m+n)^k=\dim_{\cc(q)}(V^{\otimes k})
 &\ge \sum_{\substack{\text{$T$: st. tab. of} \\
                      \quad \text{shape $\lambda$} \\
               \lambda \in H(m,n;k)}}
     \dim_{\cc(q)}(y_T(q))V^{\otimes k}) \\
 &\ge \sum_{\substack{\text{$T$: st. tab. of}
                      \quad \text{shape $\lambda$} \\
              \lambda \in H(m,n;k)}}
   \dim_{\cc}(y_T)W^{\otimes k}) \\
 & =W^{\otimes k} = (m+n)^k.
\end{split}
\end{equation}
Therefore all the equalities  in \eqref{eq:ee} should hold, and
\begin{equation*}
V^{\otimes k}=\bigoplus_{\substack{\text{$T$:
  st. tab.} \\
  \quad \text{of shape $\lambda$} \\
  \lambda \in H(m,n;k)}}
 y_T(q)V^{\otimes k}).
\end{equation*}
Now the corollary follows.
\end{proof}

NOw note that $\dim_\cc I_\lambda=\dim_{\cc(q)} I^q_\lambda=
d_\lambda^2$,  and we have
\begin{equation}\label{eq:dimeqhec}
\dim_{\cc} \bigl( \Psi(\cc S_k) \bigr)= \dim_{\cc(q)}
\bigl(\Psi_q(\hec)\bigr).
\end{equation}

Our next goal is to prove that the image $\Psi_q(\hec)$ of $\hec$
is in fact the full centralizer $\End_\uu (V^{\otimes k})$ of
$\uu$ on $V^{\otimes k}$.  The completely reducibility of
$\uu$-module $V^{\otimes k}$ and the branching rule for a tensor
product of the $\uu$ module $V$ are obtained by G.~Benkart,
S.~Kang and M.~Kashiwara using crystal graphs of $\uq(gl(m,n))$ in
\cite{bkk}.

\begin{pro}[See Proposition 3.1 in \cite{bkk}]\label{pro:reduci}
The $\uu$-module $V^{\otimes k}$ is completely reducible for all
$k\ge 1$.
\end{pro}

\begin{thm}[See Theorem 4.13 in \cite{bkk}]\label{thm:branch}
Let $\lambda_0\vdash k$ be an $(m,n)$ hook-shape. Then the tensor
product $V(\lambda_0) \otimes V(\epsilon_1)$ has the following
decomposition into irreducible $\uq(gl(m,n))$--modules:
\begin{equation} \label{eq:branch}
V(\lambda_0) \otimes V(\epsilon_1) = \bigoplus_{\lambda \in \Lambda} V(\lambda),
\end{equation}
where $\lambda$ runs over the set  $\Lambda$ of all $(m,n)$ hook-shape Young
diagrams obtained from $\lambda_o$ by adding a box to $\lambda_0$.
\end{thm}

Let $\mathfrak{g}$ be a  Lie superalgebra. A $\mathfrak{g}$-module
$V$ is \emph{irreducible} if $V$ does not have $\mathfrak{g}$
invariant $\mathbb{Z}_2$-graded subspace. Note that irreducible
modules appearing in \eqref{eq:branch} do not have any
$\mathbb{Z}_2$-graded or non-graded subspace which is
$\mathfrak{g}$ invariant.  Therefore the Schur's lemma is still
true in our case.
\begin{lem}\label{lem:schur}
Let $V(\lambda)$ and $V(\mu)$ be any two irreducible $\uu$-module
appearing in the branching rule \eqref{eq:branch}. Then
$$
\End_\uu \bigl(V(\lambda), V(\mu)\bigr)= \begin{cases}
   \cc(q) &\text{ if $\lambda=\mu$}, \\
   0  &\text{ if $\lambda \ne \mu$}.
\end{cases}
$$
\end{lem}

The branching rule for tensor products of the of $\mathfrak{U}(gl(m,n))$-module
$W=\cc^m \oplus \cc^n$ was obtained by Berele and Regev, and it is same as
\eqref{eq:branch}. The centralizer theorem
\begin{equation} \label{thm:br}
\End_{\mathfrak{U}(gl(m,n))}(W^{\otimes k}) = \Psi(\cc S_k)
\end{equation}
for the nonquantum case was also obtained by Berele and Regev \cite{br}.

Now because the branching rules for quantum and nonquantum cases
are the same, we have from Lemma~\ref{lem:schur} and
\eqref{eq:dimeqhec} that
\begin{equation}
\begin{split}
\dim_{\cc(q)}\End_{\uu}(V^{\otimes k}) &=
\dim_{\cc}\End_{\mathfrak{U}(gl(m,n))}(W^{\otimes k})  \\
&=\dim_\cc\Psi(\cc S_k)\\
&=\dim_{\cc(q)}\Psi_q(\hec).
\end{split}
\end{equation}
Thus we have the following centralizer theorem for $\uu$.
\begin{thm}\label{thm:centra}
The centralizer algebra of the actions of
$\mathfrak{U}_q(gl(m,n))$ on $V^{\otimes k}$ is the image of
Iwahori-Hecke algebra  $\hec$ under $\Psi_q$, i.e.
$$
\End_{\mathfrak{U}_q(gl(m,n))}(V^{\otimes k}) = \Psi_q(\hec).
$$
\end{thm}
Moreover the {\em double centralizer theory} gives the following:
\begin{cor}\label{cor:double}
The centralizer algebra of the action of $\hec$ on $V^{\otimes k}$
is the image of the quantized enveloping algebra $\uu$ under
$\rho^{\otimes k}$, i.e.
$$
\End_{\Psi(\hec)}(V^{\otimes k})= \rho^{\otimes k}\left(\uu\right).
$$
\end{cor}

\section{Symmetric groups and $k$-diagrams} \label{sec:little}

It is helpful to represent permutations in the symmetric group $S_k$
by  diagrams.  Consider a graph with two rows of $k$ vertices
each, one above the other, and $k$ edges such that each vertex in
the top row is incident to precisely one vertex in the  bottom row.
There is a natural one-to-one correspondence between such $k$-diagrams and
elements of the symmetric group $S_k$, which is illustrated by the following
example:
\begin{exa} \label{exa:diag}
$$
\begin{pmatrix}
1& 2& 3& 4& 5& 6 \\
3&5&6&1&4&2
\end{pmatrix}
 =\rotatebox[origin=Bc]{270}{\rotatebox[origin=tr]{90}{\scalebox{0.4}
 {\includegraphics{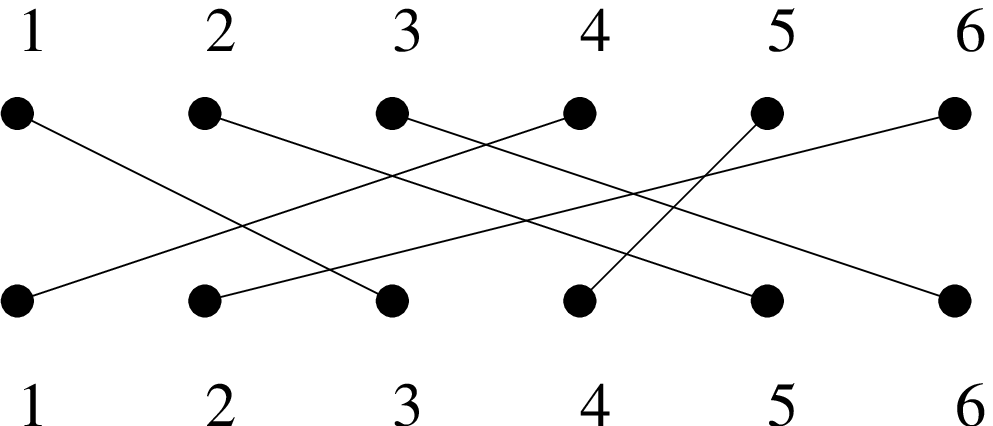}}}}.
$$
\end{exa}
Notice that the $i$th vertex in top row is incident to the
$\sigma(i)$th vertex in bottom row. We identify the generator
$s_i=(i\, i+1)$ of $S_k$ with the following diagram:
$$
s_i= \rotatebox[origin=Bc]{270}{\rotatebox[origin=tr]{90}{\scalebox{0.4}
{\includegraphics{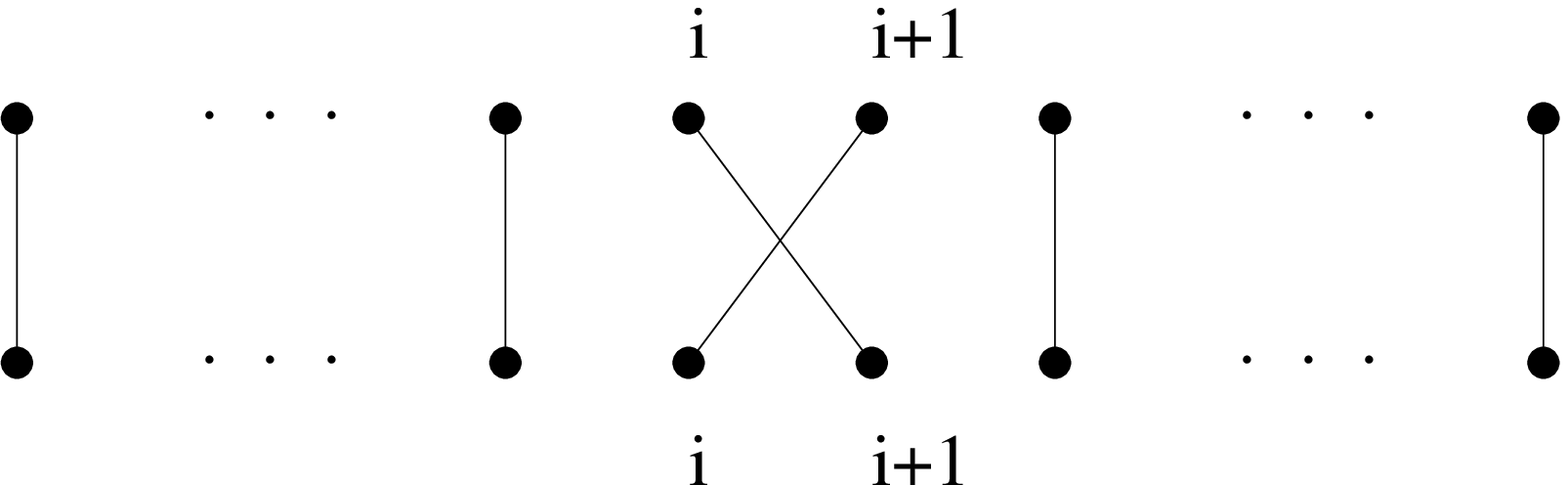}}}}.
$$

Let $d_1$ and $d_2$ be the diagrams corresponding to permutations
$\sigma_1$ and $\sigma_2$ respectively. Place $d_1$ below $d_2$
and identify the vertices in the bottom row of $d_2$ with the
corresponding vertices in the top row of $d_1$. The resulting
diagram is corresponding to the product $\sigma_1\sigma_2$.  For
example,
$$
\begin{pmatrix}
1&2&3 \\ 2&3&1
\end{pmatrix}= (12)(23) =
\rotatebox[origin=Bc]{270}{\rotatebox[origin=tr]{90}{\scalebox{0.4}
{\includegraphics{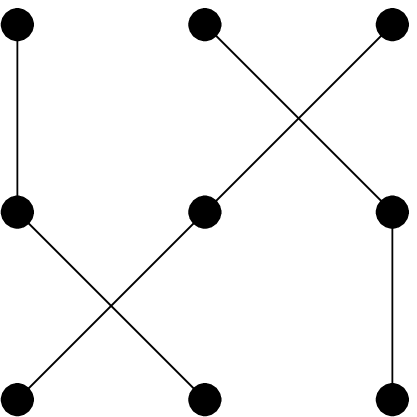}}}}
=\rotatebox[origin=Bc]{270}{\rotatebox[origin=tr]{90}{\scalebox{0.4}
{\includegraphics{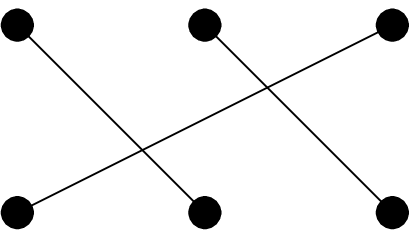}}}}.
$$
Note that we stack the left element of the product on the bottom of the
diagram and the right element on the top.

Let
\begin{equation}\label{eq:skexp}
\sigma=s_{i_1}s_{i_2}\cdots s_{i_l},
\end{equation}
be an expression of $\sigma\in S_k$. Then the {\em length}
$\ell(\sigma)$ of $\sigma \in S_k$ is the number of crossings of
edges in the $k$-diagram identified with $\sigma\in S_k$. An
expression $\sigma=s_{i_1}\cdots s_{i_j}$ of $\sigma\in S_k$ is
reduced  if $j=\ell(\sigma)$. For example, the $k$-diagram shown
in Example~\ref{exa:diag} has $9$ edge crossings, and so
$\ell(\sigma)= 9$, and
\begin{equation*}
\sigma=\begin{pmatrix}
1&2&3&4&5&6 \\ 3&5&6&1&4&2
\end{pmatrix}=s_3s_4s_2s_3s_4s_5s_1s_2s_3 ,
\end{equation*}
where the product on the right is a reduced expression for $\sigma$.

 Let $\sigma=\sigma_1 \sigma_2$ be the
product of two permutations $\sigma_1$ and $\sigma_2$.
Then $\ell(\sigma) < \ell(\sigma_1)+ \ell(\sigma_2)$ if and only if the
situation explained below using $k$-diagrams happens :
\begin{equation} \label{eq:reduction}
\rotatebox[origin=Bc]{270}{\rotatebox[origin=tr]{90}{\scalebox{0.4}
 {\includegraphics{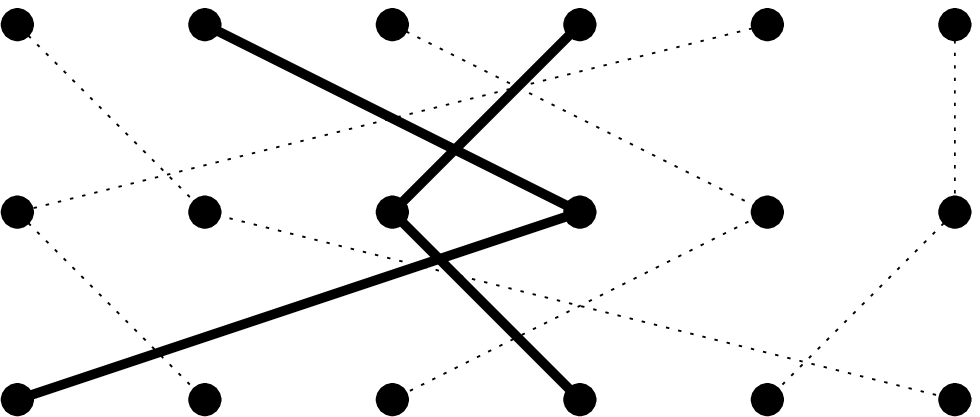}}}} \Rightarrow
\rotatebox[origin=Bc]{270}{\rotatebox[origin=tr]{90}{\scalebox{0.4}
 {\includegraphics{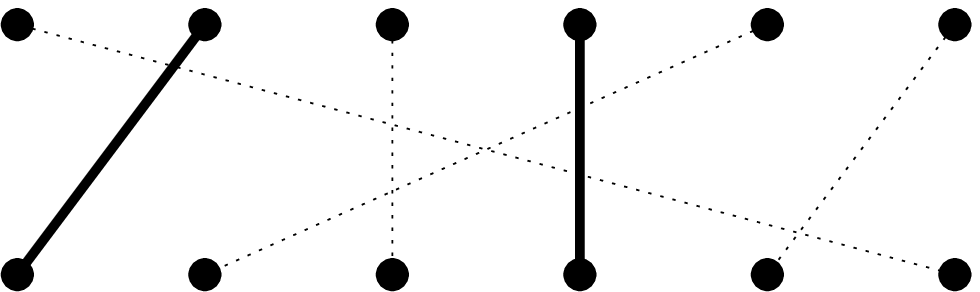}}}}.
\end{equation}
The crossings given by the darkened edges disappear in the product.

Let $\mathcal{L}$ be a set consisting of $l$ letters $\left\{a_1,
a_2,\ldots, a_l \right\}$. Let $x=x_1\cdots x_k$ be a word of
length $k$, where $x_i \in \mathcal{L}$. Then the symmetric group
$S_k$ acts on the set of all words of length $k$ by place
permutations, i.e., for $\sigma\in S_k$,
$$
\sigma(x)= x_{\sigma^{-1}(1)} \cdots x_{\sigma^{-1}(k)}.
$$
For example when $\sigma$ is the permutation in
Example~\ref{eq:reduction}, then the action of $\sigma$ on
$x=x_1\cdots x_6$ is explained as
$$
\rotatebox[origin=Bc]{270}{\rotatebox[origin=tr]{90}{\scalebox{0.4}
 {\includegraphics{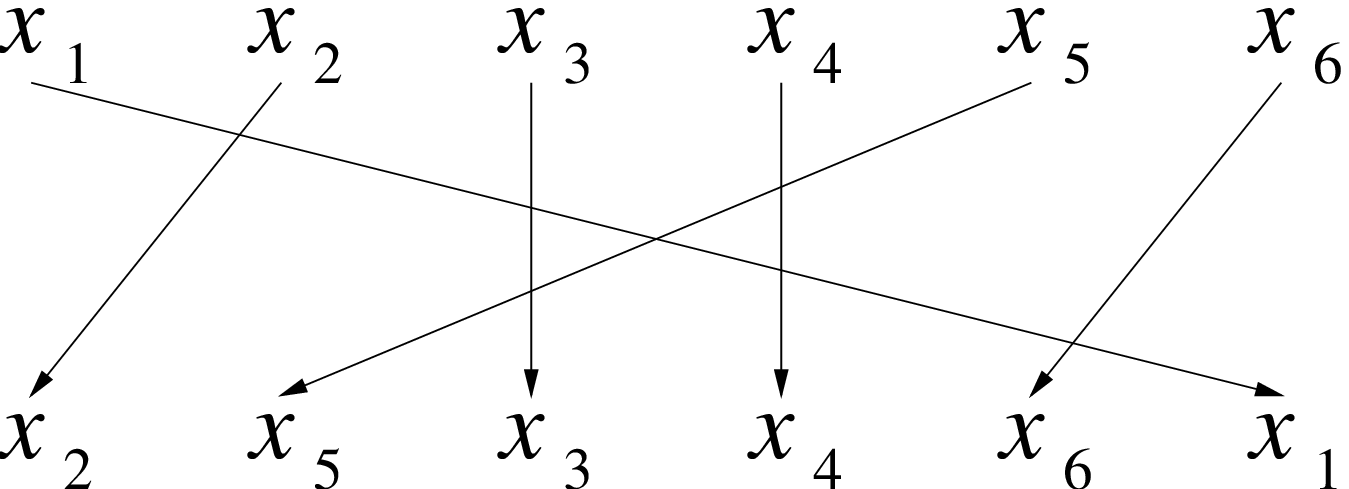}}}},
$$
so that $\sigma x= x_2 x_5 x_3 x_4 x_6 x_1$.

\section{Maximal vectors of $\mathfrak{U}_q(gl(m,n))$ modules.}
\label{sec:maximal}

In this section we construct the highest weight vectors of the
irreducible summands of $\uu$-module $V^{\otimes k}$ explicitly
using Gyoja's $q$-analogue of the Young symmetrizers.

First we note the following lemma from \cite{gyoja}.
\begin{lem}[See \cite{gyoja}.]\label{lem:trick}
Let $T$ be a standard tableau of shape $\lambda\vdash k$. Then
there exists a $\gamma \in \cc(q)$ such that
$$
e^-_\lambda \left(h(\sigma^T_-)\right)^{-1}h(\sigma^T_+)e^+_\lambda=
\gamma e^-_\lambda h(\sigma^-_+)e^+_\lambda.
$$
\end{lem}

Now let $\widehat\Pi(m,n;k)$ be
\begin{equation}\label{eq:pihat}
\widehat\Pi(m,n;k):=\left\{(\mu,\nu)\left|
\begin{array}{c}
\mu \vdash s, \nu \vdash t, \quad s+t=k, \\
\ell(\mu) \le m, \ell(\nu)\le n, \text{ and }\mu_m \ge \ell(\nu)
\end{array}
\right.\right\}.
\end{equation}
\begin{lem} [\cite{nova}]
There is a bijection between $H(m,n;k)$ and $\widehat\Pi(m,n;k)$
given  by $\lambda \mapsto (\lambda^1,\lambda^2)$, where
$$
\lambda^1=(\lambda_1,\ldots, \lambda_m),\quad \text{ and }
\lambda^2=(\lambda^2_1,\ldots, \lambda^2_n),
$$
such that $\lambda^2_j=\max\{\lambda^*_j-m,0\}$, for
$j=1,\ldots,n$.
\end{lem}

For a standard tableau $T$ of shape $\lambda=(\lambda^1,\lambda^2) \in
H(m,n;k)$, we let $T_{\lambda^1}$ be the subtableau of $T$ of
shape $\lambda^1$ and let $T_{\lambda^2}$ be the conjugate of the
skew tableau $T/\ T_{\lambda^1}$. Then we associate to $T$ a simple
tensor $w_T=v_1\otimes \cdots \otimes v_k$ in $V^{\otimes k}$ which is
defined by
\begin{equation*}
v_l = \begin{cases}
t_i \quad &\text{ if $l$ is in the $i$th row of $T_{\lambda^1}$}, \\
u_j \quad &\text{ if $l$ is in the $j$th row of $T_{\lambda^2}$}.
\end{cases}
\end{equation*}
Note that the weight of $w_T$ is
\begin{equation*}
\lambda=\lambda_1 \epsilon_1 + \cdots \lambda_{m}
\epsilon_{m}+\lambda^2_1\delta_1+\cdots +\lambda^2_n\delta_n \in
\Gamma_+ ,
\end{equation*}
where $\Gamma^+$ is defined in \eqref{eq:polwt}. For a partition
$\lambda \in H(m,n;k)$, we denote $w_\lambda^+:= w_{S^+_\lambda}$
and $w_\lambda^-:=w_{S^-_\lambda}$.

\begin{exa}
Suppose $m=2$ and $n=3$. Let $\lambda = (4,2,2,1,1)$. Then
$\lambda \in H(2,3;10)$.
$$
\lambda=
\rotatebox[origin=Bc]{270}{\rotatebox[origin=tr]{90}{\scalebox{0.8}
{\includegraphics{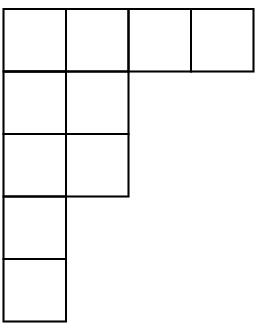}}}}
$$
The corresponding pair $(\lambda^1,\lambda^2)$ is given by
$\lambda^1=(4,2)$ and $\lambda^2=(3,1)$ so that
$$
\lambda \mapsto
\left( \;
\rotatebox[origin=Bc]{270}{\rotatebox[origin=tr]{90}{\scalebox{0.8}
{\includegraphics{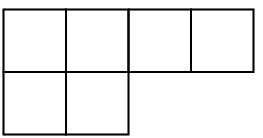}}}}\; , \;
\rotatebox[origin=Bc]{270}{\rotatebox[origin=tr]{90}{\scalebox{0.8}
{\includegraphics{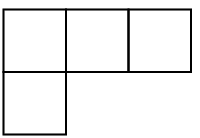}}}}
\; \right)
$$
Let $T$ be a standard tableau of shape $\lambda$ such that
$$
T = \rotatebox[origin=Bc]{270}{\rotatebox[origin=tr]{90}{\scalebox{0.8}
{\includegraphics{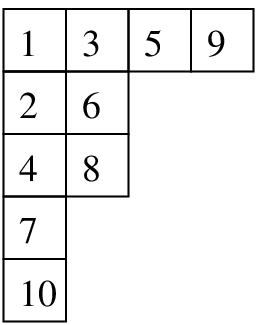}}}}
$$
Then
$$
T_{\lambda^1} =
\rotatebox[origin=Bc]{270}{\rotatebox[origin=tr]{90}{\scalebox{0.8}
{\includegraphics{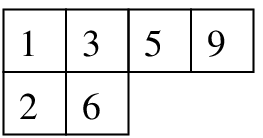}}}}\, , \qquad
T_{\lambda^2} =
\rotatebox[origin=Bc]{270}{\rotatebox[origin=tr]{90}{\scalebox{0.8}
{\includegraphics{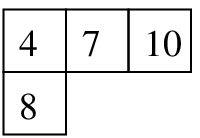}}}}\, .
$$
Then the simple tensor $w_T\in V^{\otimes 10}$ is given as
$$
w_T = t_1\otimes t_2\otimes t_1\otimes u_1\otimes t_1\otimes
t_2 \otimes u_1 \otimes u_2 \otimes t_1 \otimes u_1.
$$
\end{exa}

For a simple tensor $w=v_1\otimes \cdots\otimes v_k$, by word $w$
we mean the word $i_1i_2\cdots i_k$, where $i_j=l$ if $v_j=t_l$
and $i_j=\overline l$ if $v_j=u_l$. As we noted in
Section~\ref{sec:little}, a permutation $\sigma \in S_k$ acts on
words in $v\in V^{\otimes k}$ by place permutations .

\begin{thm} \label{thm:final}
Let $\lambda$ be a partition in $H(m,n;k)$ and $T$ be a standard
tableau of shape $\lambda=(\lambda_1,\lambda_2)$. Let $v_+=
y_T(q)h(\sigma^T_+)w_\lambda^+$. Then $v_+$ is a maximal vector in
$y_T(q)(V^{\otimes k})$ of weight $\lambda=\lambda_1\epsilon_1 +
\cdots + \lambda_m\epsilon_m+ \lambda^2_1\delta_1 + \cdots +
\lambda^2_n\delta_n \in \Gamma^+$. Hence $y_T(q)(V^{\otimes k})
\cong V(\lambda)$, the irreducible $\uu$-module with highest
weight $\lambda$.
\end{thm}
\begin{proof}
First note that when $q$ goes to $1$, $v_+=y_T(q)h(\sigma^T_+)
w_\lambda^+$ goes to a maximal vector $y_T \sigma^T_+w_\lambda^+$
of $y_T (W^{\otimes k})$ in the classical case (see \cite{nova}).
Therefore we know that $v_+$ is a nonzero vector.

Next observe that the weight of $y_T(q) h(\sigma^T_+)w_\lambda^+$ is
same as the weight of $w_\lambda^+$, which is $\lambda=
\lambda_1\epsilon_1 +\cdots +\lambda_m\epsilon_m+\lambda^2_1
\delta_1+ \cdots + \lambda^2_n\delta_n$ by construction, where
$\lambda^2_j=\max\{\lambda^*_j-m,0\}$.

Now let's prove that $v_+$ is annihilated by the action of $\rho^{\otimes k}(E_i)$
for $i=1,\ldots, m+n-1$. Note
\begin{align}\label{eq:bun}
v_+=&y_T(q)h(\sigma^T_+)w_\lambda^+ \\
=&\frac{1}{\xi}
h(\sigma^T_-)e^-_\lambda\bigl(h(\sigma^T_-)\bigr)^{-1}h(\sigma^T_+)e^+_\lambda
\bigl(h(\sigma^T_+)\bigr)^{-1}h(\sigma^T_+)w_\lambda^+ \notag\\
=&\frac{1}{\xi}
h(\sigma^T_-)e^-_\lambda\bigl(h(\sigma^T_-)\bigr)^{-1}h(\sigma^T_+)e^+_\lambda
w_\lambda^+ \notag.
\end{align}

{\em (i) First we consider the case $1\le i \le m$.}

In this case $\rho^{\otimes k}(E_i)$ maps $t_{i+1} \mapsto t_i$
for $i=1,\ldots, m$ or $u_1 \mapsto t_m$, because $\rho^{\otimes
k} (E_i)= E_{i,i+1}$. We fix a special standard tableau
$\widetilde T$ corresponding to partition $\lambda\in H(m,n)$
which is the shape of $T$. We give an example below and do not
bother to give the precise definition of $\widetilde T$.  If
$$
\lambda=
 \rotatebox[origin=Bc]{270}{\rotatebox[origin=tr]{90}{\scalebox{0.4}
 {\includegraphics{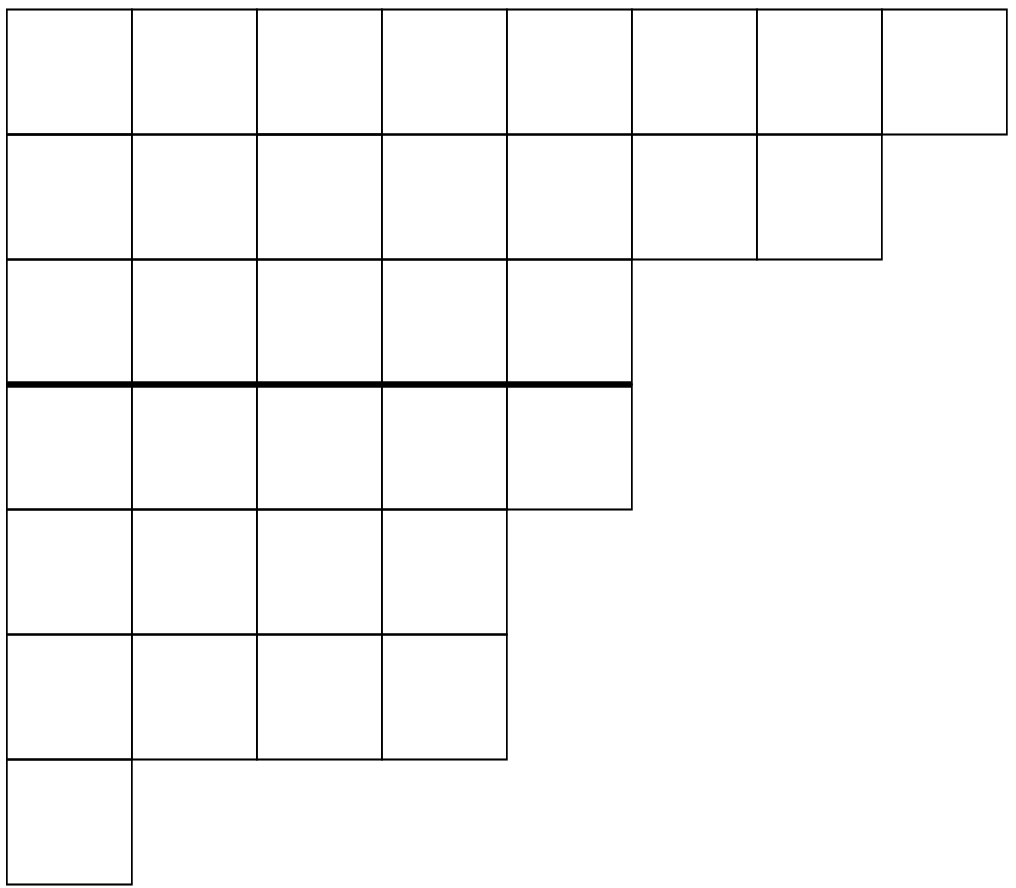}}}}\in H(3,5),
$$
then
$$
\widetilde T
 =\rotatebox[origin=Bc]{270}{\rotatebox[origin=tr]{90}{\scalebox{0.4}
 {\includegraphics{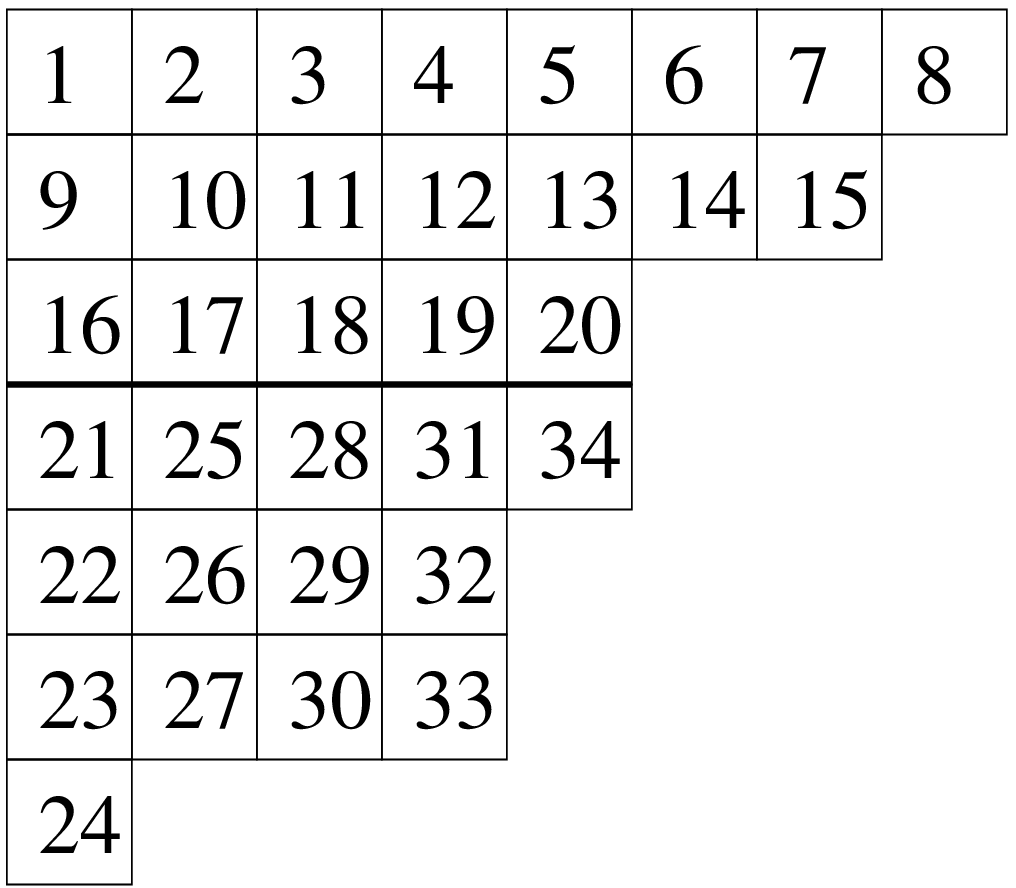}}}}.
$$
Note that the entries in $\widetilde T$ increase across rows
from left to right for the first $m(=3)$ rows, then entries increase
down the columns.

{F}rom Lemma~\ref{lem:trick} with $T=\widetilde T$, we have that
\begin{equation} \label{eq:ttt}
e^-_\lambda\left(h(\sigma^{\widetilde T}_-)\right)^{-1}
h(\sigma^{\widetilde T}_+)
e^+_\lambda=c_1 e^-_\lambda h(\sigma^-_+)e^+_\lambda ,
\end{equation}
for some $c_1\in \cc(q)$. Note we know $c_1$ is nonzero because the left
side of \eqref{eq:ttt} is nonzero by specializing $q\mapsto 1$.
Combining Lemma~\ref{lem:trick} and \eqref{eq:ttt} we have
\begin{equation}\label{eq:res}
e^-_\lambda \bigl(h(\sigma^T_-)\bigr)^{-1}h(\sigma^T_+)e^+_\lambda = c
e^-_\lambda \left(h(\sigma^{\widetilde T}_-)\right)^{-1}h(\sigma^{\widetilde T}_+)
e^+_\lambda ,
\end{equation}
for some nonzero $c\in \cc(q)$.
Now  from \eqref{eq:bun} and \eqref{eq:res}
\begin{align*}
v_+=& \frac{1}{\xi} h(\sigma^T_-)e^-_\lambda \bigl(h(\sigma^T_-)\bigr)^{-1}
h(\sigma^T_+)e^+_\lambda  w_\lambda^+\\
=& c'(q) h(\sigma^T_-)\underbrace{
e^-_\lambda \bigl(h(\sigma^{\widetilde T}_-)\bigr)^{-1}
h(\sigma^{\widetilde T}_+)e^+_\lambda  w_\lambda^+}_{(\ast)},
\end{align*}
for some $c'(q)\in \cc(q)$.
We will show that $(\ast)$ is a linear sum of simple tensors which are
killed by the actions of $E_{i\, i+1}$ for $i\le m$.

Write $I_0=\{ 1,\ldots, m\}$ and $I_1=\{\overline 1, \ldots,
\overline n\}$. Then we may write the simple tensor $w_\lambda^+$ as
 \begin{equation}\label{eq:1}
 \begin{split}
w_\lambda^+=& \underbrace{1\otimes \cdots \otimes
1}_{\lambda_1}\otimes \underbrace{2\otimes \cdots \otimes
2}_{\lambda_2}\otimes  \cdots \otimes \underbrace{m\otimes \cdots
\otimes m }_{\lambda_m} \\
 & \quad\otimes  \overline 1\otimes \overline 2\otimes \cdots\otimes
 \overline\lambda_{m+1} \quad \otimes \overline 1 \otimes
 \cdots \otimes \overline \lambda_{m+2}\otimes \cdots\otimes  \overline 1
 \otimes \cdots \otimes \overline \otimes \lambda_{\ell(\lambda)}.
 \end{split}
 \end{equation}

Recall that $\displaystyle e^+_\lambda = \sum_{\sigma \in R(S_+)}
h(\sigma)$.  Let's fix a $\sigma\in R(S_+)$. Then
$h(\sigma)w_\lambda^+$ is a linear sum of simple tensors
$w'_{\lambda_+}$ where $w'_{\lambda_+}$ are the same as
$w_\lambda^+$ up to scalar multiplications except the orders of
the entries $\overline i$ in $I_1=\{ \overline 1, \ldots,\overline
n\}$ are changed. This is because if $\sigma$ moves entries in the
$j$th row of $S_+$, $j\le m$, then the action of $h(\sigma)$ on
$w_\lambda^+$ does not create any new simple tensors because $\hr$
maps $t_i\otimes t_i$ to $q^2t_i\otimes t_i$.  And if $\sigma$
moves entries in the $j$th row of $S_+$, $j\ge m+1$, then
$h(\sigma)w_\lambda^+$ is a linear sum of $w'_{\lambda_+}$ where
$w'_{\lambda_+}$ are simple tensors same as $w_\lambda^+$ except
the orders of entries in $I_1$ are changed because
 $$
 \hr(u_i\otimes u_j) =
 \begin{cases} -q u_j\otimes u_i + (q^2-1)u_i\otimes u_j
         &\text{ if $i<j$}, \\ -q u_j\otimes u_i
         &\text{ if $i>j$},
 \end{cases}
  $$
where the $u_i$'s are basis vectors of $V_{\bar 1}$. Thus
$h(\sigma)$ maps \eqref{eq:1} to a linear sum of simple tensors
 \begin{equation}\label{eq:2}
\underbrace{1\otimes \cdots \otimes 1}_{\lambda_1}\otimes
\underbrace{2\otimes \cdots\otimes 2}_{\lambda_2} \otimes \cdots\otimes
\underbrace{m\otimes \cdots\otimes  m }_{\lambda_m}
\otimes \overline *\otimes \overline *\otimes \cdots\cdots
\otimes \overline * .
 \end{equation}

Now note that $\sigma^{\widetilde T}_+$ is the permutation which
transforms $S_ +$ to $\widetilde T$, and $\sigma^{\widetilde T}_+$
maps the word $w_\lambda^+$ to the word $w_{\widetilde T}$ by
place permutation. Moreover, the permutation $\sigma^{\widetilde
T}_+$ does not move entries from $I_0=\{1,\ldots,m\}$. Thus the
action of the Hecke element $h(\sigma^{\widetilde T}_+)$ on
\eqref{eq:2} produce a linear sum of simple tensors
$w'_{\widetilde T}$ where $w'_{\widetilde T}$ are the same as
$w_{\widetilde T}$ except only the orders of entries in $I_1$ are
changed, which are in fact same as $w^+_\lambda$ except only the
orders of entries in $I_1$ are changed. Thus we may also write
$w'_{\widetilde T}$ just like \eqref{eq:2}.
 \begin{equation} \label{eq:4}
\underbrace{1\otimes \cdots \otimes 1}_{\lambda_1}\otimes
\underbrace{2\otimes \cdots\otimes 2}_{\lambda_2} \otimes \cdots\otimes
\underbrace{m\otimes \cdots\otimes  m }_{\lambda_m}
\otimes \overline *\otimes \overline *\otimes \cdots\cdots
\otimes \overline * .
 \end{equation}


Our next goal is to show $\left(h(\sigma^{\widetilde T}_{-})
\right)^{-1}$ maps \eqref{eq:4} to a scalar multiple of
$w'_{\lambda_-}$, where $w'_{\lambda_-}$ is a simple tensor which
is the same as $w_\lambda^-$ except the order on entries
$\overline i$'s is different. Write the simple tensor
$w_\lambda^-$ as
\begin{equation}\label{eq:3}
\begin{split}
w_\lambda^-=& 1 \otimes  2\otimes \cdots \otimes m \otimes
{\underbrace{\overline 1\otimes \cdots\otimes  \overline
1}_{(\lambda^2)_1}}\otimes 1\otimes \cdots \otimes m \otimes
{\underbrace{\overline 2\otimes \cdots \otimes \overline
2}_{(\lambda^2)_2}} \\
&\quad \otimes  \cdots\otimes  1 \otimes \cdots \otimes \lambda_2
\otimes \underbrace{1\otimes \cdots \otimes
1}_{\lambda_1-\lambda_2}
\end{split}.
\end{equation}

Now note that
$\sigma^{-}_{\widetilde T}$
maps the word $w_{\widetilde T}$ to the word $w_\lambda^-$ by place
permutation. We decompose $\sigma^{-}_{\widetilde T}=
\sigma_{\lambda_1}\cdots \sigma_1$, so that $\ell\left(
\sigma^{-}_{\widetilde T}\right)= \ell\left(\sigma_{\lambda_1}\right)
+ \cdots +\ell\left(\sigma_1\right)$ as in the followings ways:

First we define a sequence $\widetilde T_0=\widetilde T,
\widetilde T_1,\ldots, \widetilde T_{\lambda_2-1}, \widetilde
T_{\lambda_2} = S_\lambda^-$ of standard tableux  such that
$\widetilde T_i$ is a  standard tableau of shape $\lambda $ whose
entries increase by one down the first $i$ columns, and then other
entries increase by one just like they are in the $\widetilde{T}$
for the rest of columns, i.e. entries increase by one across  rows
for first $m$ rows then entries increase by on down the columns.
Then we define $\sigma_i := \sigma_{\widetilde
T_{i-1}}^{\widetilde T_i}$.

We give an example to explain our idea.  If
\allowdisplaybreaks
\begin{equation*}
\lambda=(\lambda^1,\lambda^2)
=\rotatebox[origin=Bc]{270}{\rotatebox[origin=tr]{90}
{\scalebox{0.35}{\includegraphics{foo2.eps}}}} \in H(3,5),
\end{equation*}
then $\widetilde T_0=\widetilde T$, $\widetilde T_1$,
$\widetilde T_2$, $\widetilde T_3$, $\widetilde T_4$,
$\widetilde T_5$, $\widetilde T_6$, and
$\widetilde T_7=S^-_\lambda$ are
\begin{alignat*}{2}
&\widetilde T_0
=\rotatebox[origin=Bc]{270}{\rotatebox[origin=tr]{90}
{\scalebox{0.35}{\includegraphics{tildh38.eps}}}},
\qquad
&&\widetilde T_1
=\rotatebox[origin=Bc]{270}{\rotatebox[origin=tr]{90}
{\scalebox{0.35}{\includegraphics{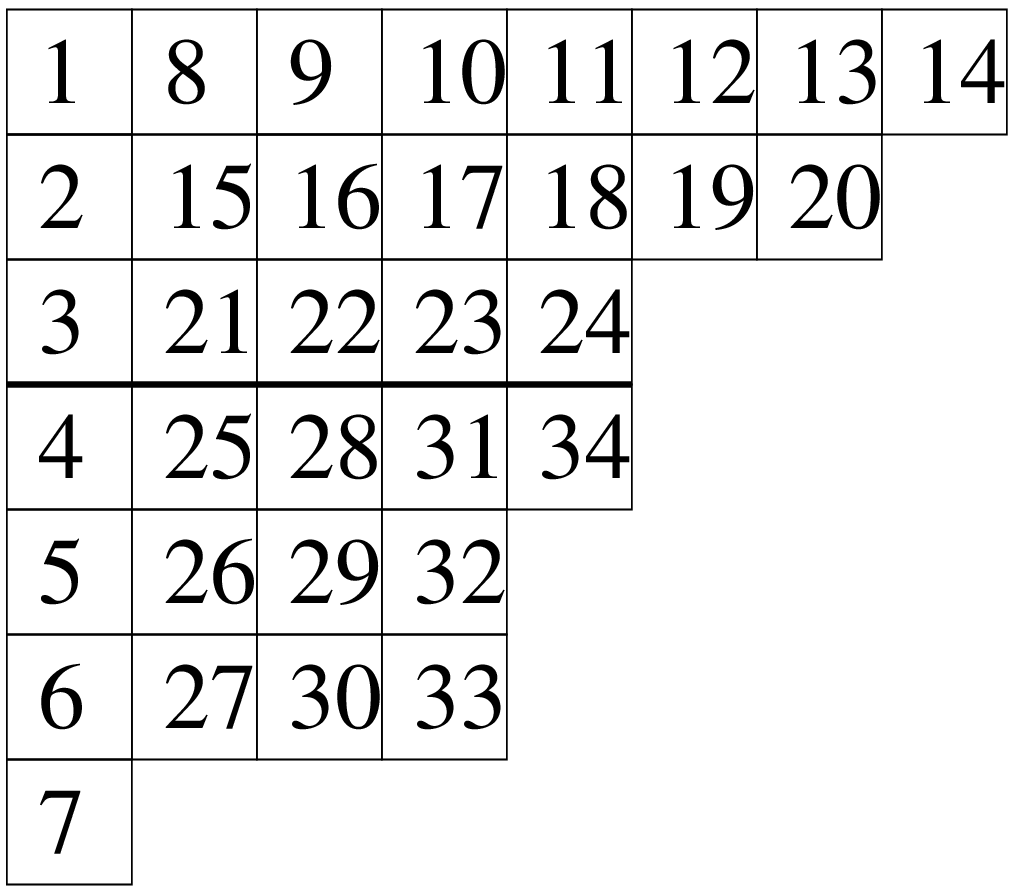}}}} \\
& &&\\
&\widetilde T_2
=\rotatebox[origin=Bc]{270}{\rotatebox[origin=tr]{90}
{\scalebox{0.35}{\includegraphics{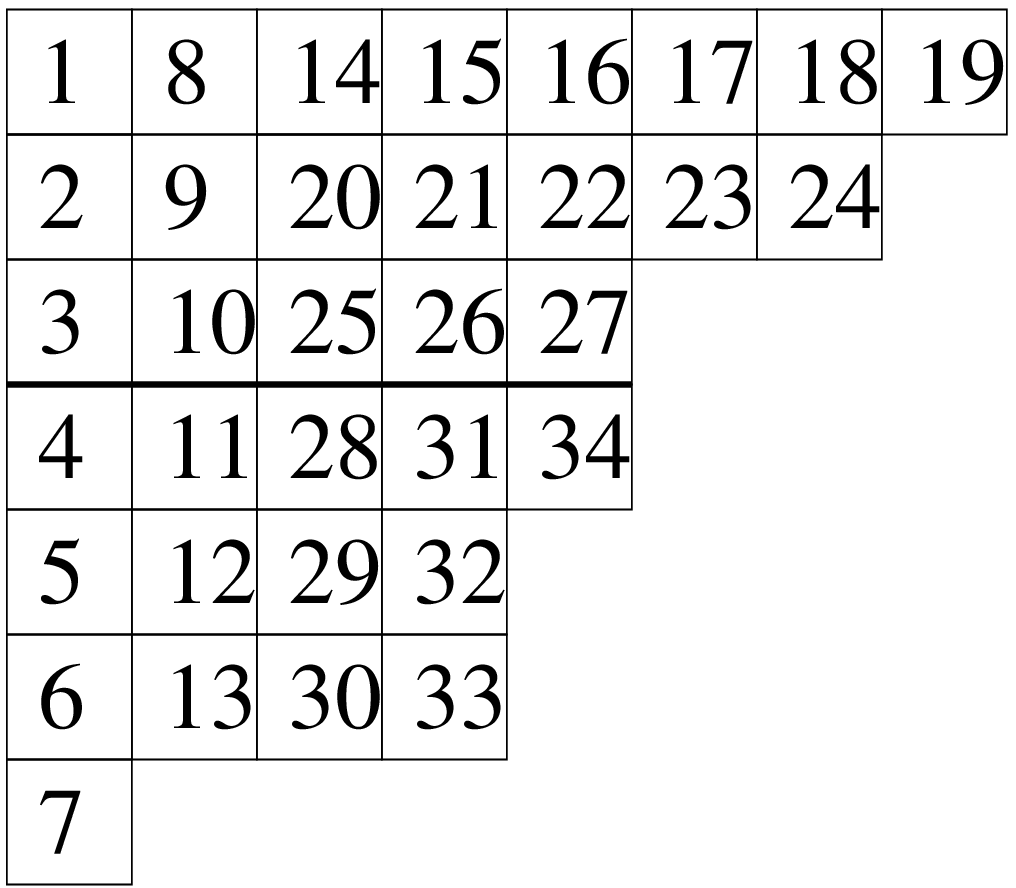}}}},
&&\widetilde T_3
=\rotatebox[origin=Bc]{270}{\rotatebox[origin=tr]{90}
{\scalebox{0.35}{\includegraphics {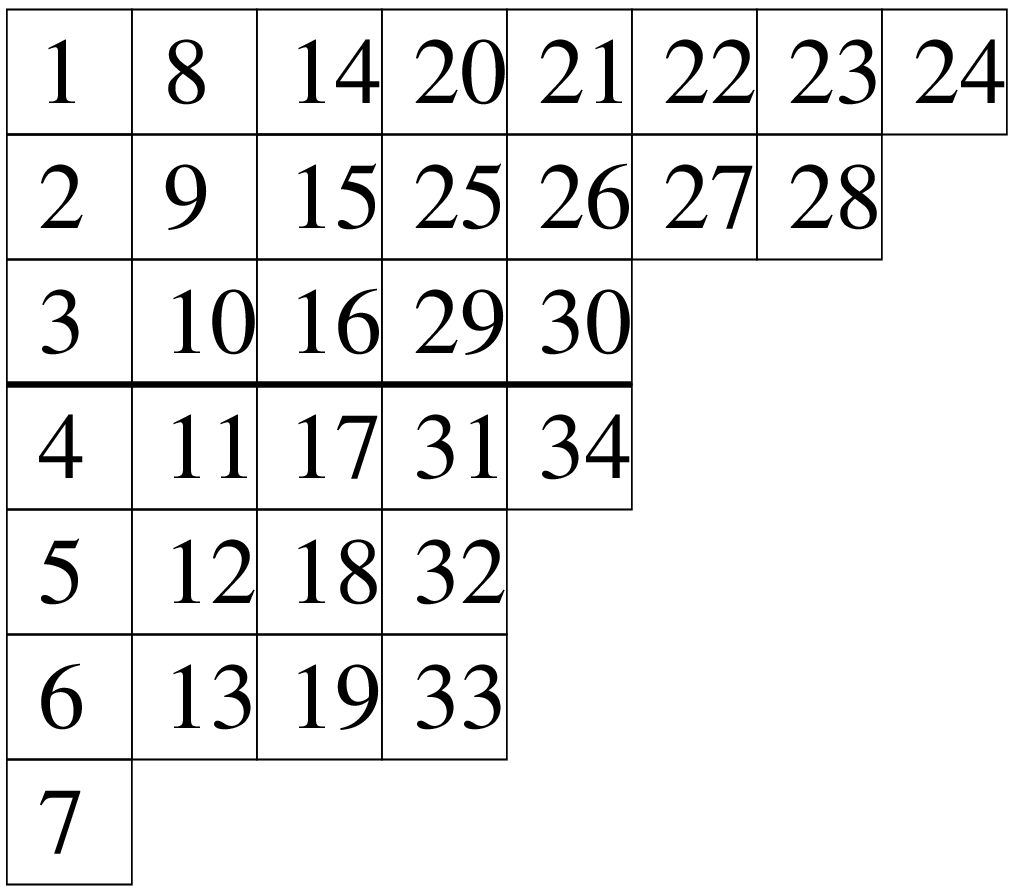}}}}, \\
& &&\\
& &\vdots &\\
& &&\\
&\widetilde T_6
=\rotatebox[origin=Bc]{270}{\rotatebox[origin=tr]{90}
{\scalebox{0.35}{\includegraphics{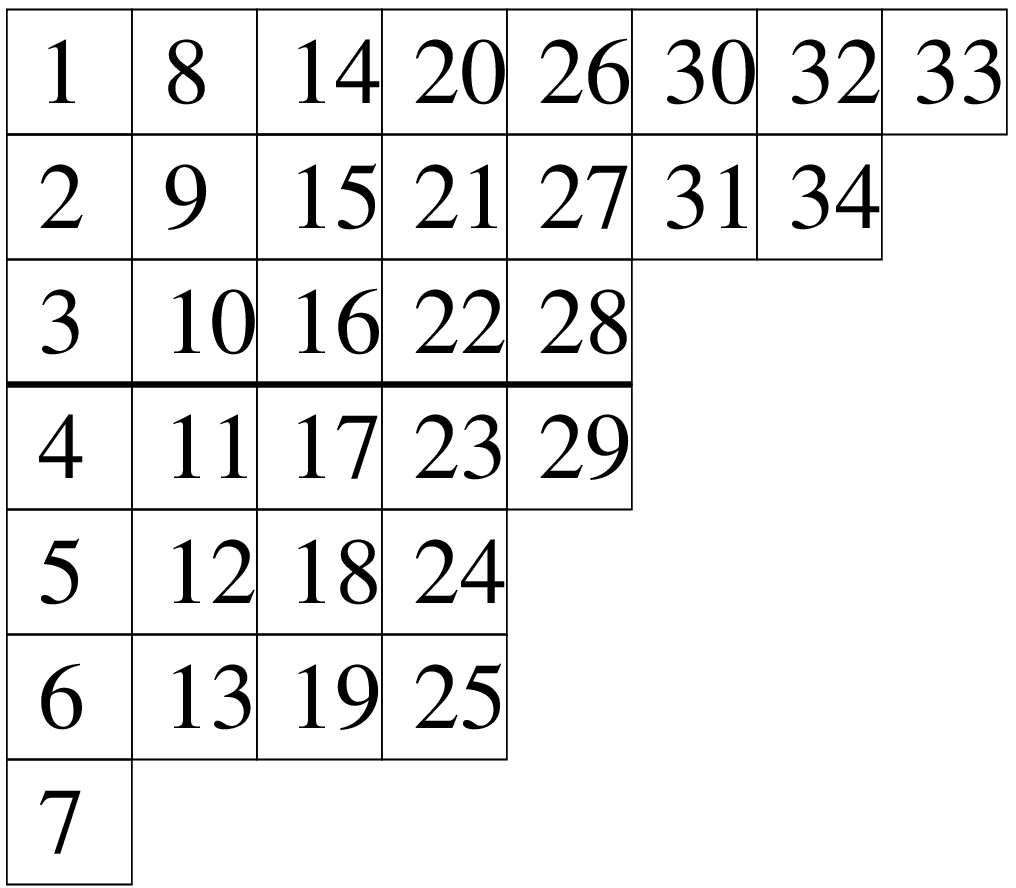}}}},
&&\widetilde T_7
=\rotatebox[origin=Bc]{270}{\rotatebox[origin=tr]{90}
{\scalebox{0.35}{\includegraphics{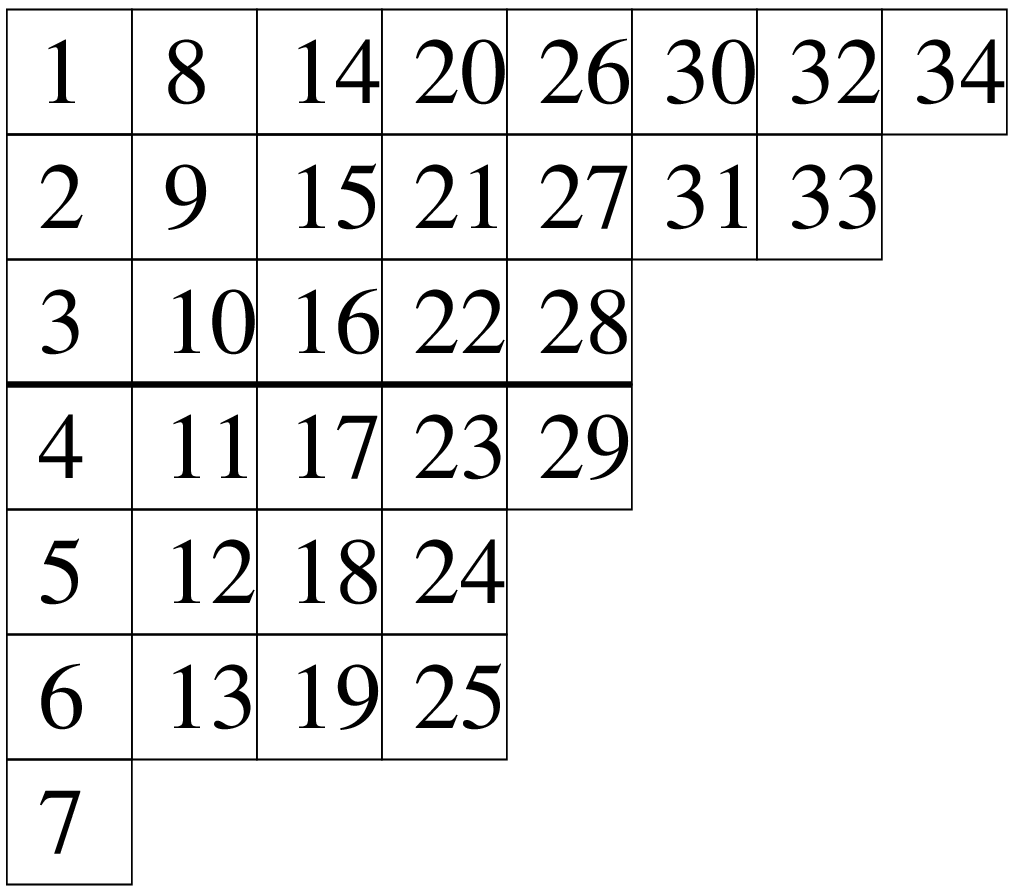}}}}.
\end{alignat*}

And our decomposition $\sigma^{-}_{\widetilde T}= \sigma_7 \cdots
\sigma_1$ is illustrated in Figure~\ref{fig:one}.

  \begin{figure}[!ht]
  \begin{center}
  \scalebox{0.27}{\includegraphics{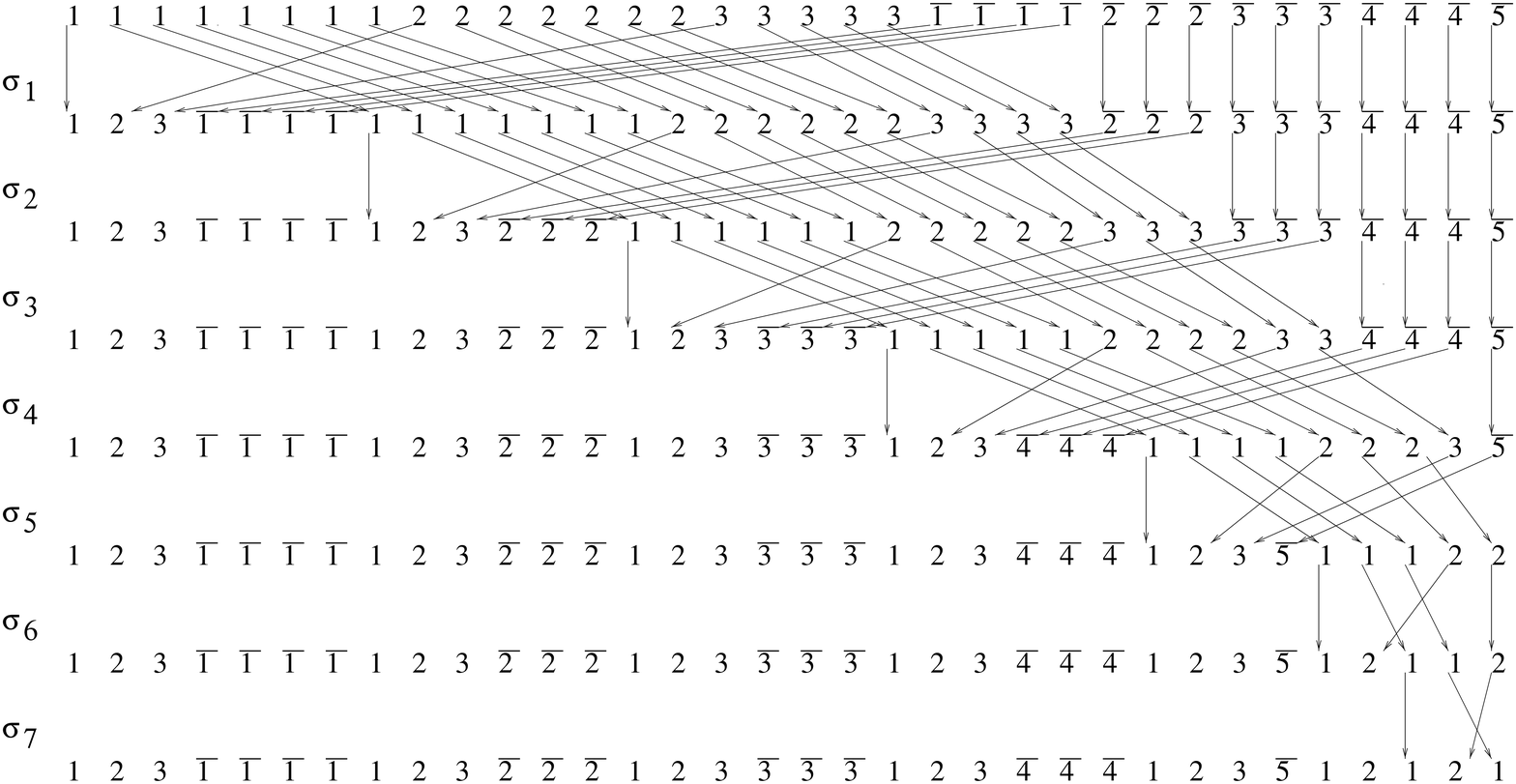}}
  \caption{$\sigma^{-}_{\widetilde T}= \sigma_7 \cdots
\sigma_1$} \label{fig:one}
   \end{center}
   \end{figure}

Note the $\overline i$'s entries of our $w'_{\widetilde T}$ and
$w'_{\lambda_-}$ are not nicely ordered as in Figur~\ref{fig:one},
but we don't need to consider those orders.

Now we notice that the situation explained in \eqref{eq:reduction}
does not happen in decomposition $\sigma^{-}_{\widetilde T}=
\sigma_7 \cdots \sigma_1$, because once an entry moves to the
left, then it is fixed by the following series of transformations
$\sigma_l$ so that it goes straight down, and it does not produce
any further crossings of edges. This is clear from
Figure~\ref{fig:one}. Thus
\begin{equation*}
\ell(\sigma_{\lambda_2}\cdots \sigma_1) =\ell(\sigma_{\lambda_2})
+ \cdots + \ell(\sigma_1),
\end{equation*}
and so,
\begin{equation*}
h(\sigma_{\lambda_2}\cdots \sigma_1) =h(\sigma_{\lambda_2})
\cdots h(\sigma_1).
\end{equation*}

Next we decompose each $\sigma_l$ into a product of
transpositions. For example a decomposition for $\sigma_4$ in
Figure~\ref{fig:one} is explained in Figure~\ref{fig:two}.

  \begin{figure}[!ht]
\begin{center}
\rotatebox[origin=Bc]{270}{\rotatebox[origin=tr]{90}
{\scalebox{0.35}{\includegraphics{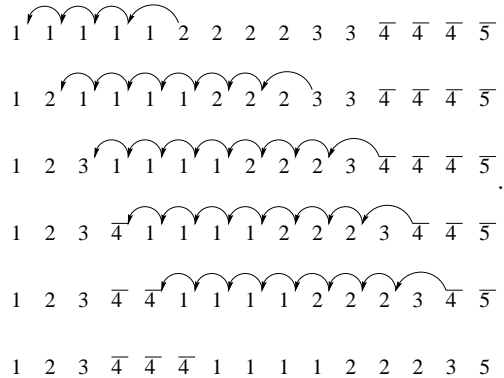}}}}.
\caption{Decomposition of $\sigma_4$ into a product of
transpositions} \label{fig:two}
\end{center}
  \end{figure}

If each $\sigma_l$ is expressed as a product
$\tau_{s_1}\tau_{s_2}\cdots \tau_{s_l}$ of transpositions as shown
in Figure~\ref{fig:two}, then the expression is reduced for the
same reason as decomposition in Figure~\ref{fig:one} is reduced.
So Hecke element $h(\sigma_l)$ is also a product of $h_s$'s which
correspond to the transpositions in the decomposition. Note also
that we only exchange $i$ and $i-1$ or $i$ and $\overline j$
during the process of applying the above place permutations.

Because $\hr$ maps $t_j\otimes t_i$ to $q t_i\otimes t_j$ for $i<j$
and $u_j\otimes t_i$ to $q t_i\otimes u_j$, we have
\begin{alignat*}{3}
\hr^{-1}: \quad &t_i\otimes t_j &&\mapsto q^{-1}t_j\otimes t_i,
 &\quad i<j&, \\
&t_i\otimes u_j &&\mapsto q^{-1}u_j\otimes t_i, \qquad &\text{ for all
$i,j$}&.
\end{alignat*}
Thus the actions of $(h_s)^{-1}$ coming from the decomposition of
$\left(h(\sigma^{\widetilde T}_-)\right)^{-1}$ on the simple
tensor $w'_{\widetilde T}$ are the same as the place permutations
except for scalar multiples.  Thus we have
$\left(h(\sigma^{\widetilde T}_-)\right)^{-1}$ maps \eqref{eq:4}
to $w'_{\lambda_-}$ except for a scalar multiple of a power of
$q^{-1}$. Hence $\bigl(h(\sigma^{\widetilde T}_-)\bigr)^{-1}
h(\sigma^{\widetilde T}_+)e^+_\lambda w_\lambda^+$ is a linear sum
of
\begin{equation}\label{eq:5}
\begin{split}
& 1\otimes  2\otimes\cdots\otimes m\otimes {\underbrace{\overline
* \otimes\cdots \overline *}_{(\lambda^2)_1}}\otimes
1\otimes\cdots \otimes m \otimes {\underbrace{\overline
*\otimes\cdots\otimes
\overline *}_{(\lambda^2)_2}} \\
&\quad \otimes\cdots\otimes 1\otimes \cdots\otimes \lambda_2 \otimes
\underbrace{1\otimes\cdots\otimes 1}_{\lambda_1-\lambda_2}
\end{split}
\end{equation}

Now we show that $e^-_\lambda w'_{\lambda_-}$ is killed by the
action of $\rho^{\otimes k}\left(E_i\right)$ for $i\le m$. The
action of $E_i$ on $V^{\otimes k}$ commute with the actions of the
Hecke algebra, so that $\rho^{\otimes k}(E_i) e_\lambda^-
w'_{\lambda_-} = e_\lambda^- \rho^{\otimes k}(E_i)
w'_{\lambda_-}$. Note $\rho^{\otimes k} (E_i) w'_{\lambda_-}$ is a
linear sum of simple tensors $\theta_\alpha$ such that each
$\theta_\alpha$ has a tensor factor where $t_i$ has been changed
to $t_{i-1}$ or $u_1$ has been changed to $t_m$. If the first case
happens, then there is a $(j\, j+1)\in C(S_-)$ such that $h_j
\theta_\alpha=q^2\theta_\alpha$ (note $\hr (t_i\otimes t_i)=
q^2t_i\otimes t_i$). So
 \begin{align*}
e^-_\lambda \theta_\alpha=&q^{-2}e^-_\lambda h_j \theta_\alpha \\
=&q^{-2}(-1)e^-_\lambda \theta_\alpha.
 \end{align*}
Hence we obtain $e^-_\lambda \theta_\alpha=0$ as expected.

If the second case happens, then for some $1\le a,b \le k$, where
$a,b$ are in the same column of $S_-$, the vectors in the $a$'th
and $b$'th tensor slots are both $t_m$. Note the vectors between
the $a$'th and $b$'th tensor slots are from $\{u_1,\ldots, u_n\}$.
Consider Hecke element $h_ah_{a+1}\cdots h_{b-2}h_{b-1}
(h_{b-2})^{-1} \cdots (h_a)^{-1} \in \hec$. Without loss of
generality we assume $a=1$. Then
\begin{align*}
h_ah_{a+1}\cdots& h_{b-2}h_{b-1} (h_{b-2})^{-1} \cdots (h_a )^{-1}
(t_m\otimes u_i \otimes \cdots \otimes u_j \otimes t_m)\\
&= q^{-(b-a-2)}h_ah_{a+1}\cdots h_{b-2}h_{b-1} (u_i\otimes \cdots
\otimes u_j \otimes t_m \otimes t_m ) \\
&= q^{-(b-a-2)}q^2 h_ah_{a+1}\cdots h_{b-2} (u_i\otimes \cdots \otimes
u_j \otimes t_m \otimes t_m ) \\
&= q^{-(b-a-2)}q^2 q^{(b-a-2)} (t_m\otimes u_i\otimes \cdots \otimes
u_j \otimes t_m ) \\
&=q^2 (t_m \otimes u_i \otimes\cdots \otimes  u_j\otimes t_m ) ,
\end{align*}
which explains that
$$
h_ah_{a+1}\cdots h_{b-2}h_{b-1} (h_{b-2})^{-1} \cdots (h_a )^{-1}
\theta_\alpha = q^2\theta_\alpha.
$$
Note $e^-_\lambda h_j=- e^-_\lambda $ and $e^-_\lambda (h_j)^{-1}=-
e^-_\lambda $ if $s_j\in C(S_-)$, so that
\begin{align*}
e^-_\lambda \theta_\alpha &=
q^{-2} e^-_\lambda  h_ah_{a+1}\cdots h_{b-2}h_{b-1} (h_{b-2})^{-1}
\cdots (h_a )^{-1} \theta_\alpha \\
&= -q^{-2} e^-_{\lambda}\theta_\alpha.
\end{align*}
Therefore we have $e^-_\lambda \theta_\alpha=0$ again this case.
Hence $E_{i \, i+1} e_\lambda^- w'_\lambda = 0$. And we have shown
here $\rho^{\otimes k}(E_i)v_+=0$ for $1\le i\le m$.

\medskip

{\em (ii) Second we consider the case $m<i\le  m+n-1$.}

This case is somewhat easier than the other case. First note
$E_{i, i+1}w_\lambda^+$ is a linear sum of simple tensors
$\theta_\alpha$ such that one of $u_{i+1}$ in tensor slots of
$w_\lambda^+$ is changed to $u_i$. Then for some $(j\, j+1)\in
R(S_+)$, $h_j \theta_\alpha= - \theta_\alpha$ because $\hr
(u_i\otimes u_i)=-u_i \otimes u_i$. Thus
\begin{align*}
e^+_\lambda  \theta_\alpha&=- e^+_\lambda  h_j\theta_\alpha  \\
&=-q^2 e^+_\lambda  \theta_\alpha.
\end{align*}
Thus we have $e^+_\lambda \theta_\alpha=0$ and we have
$\rho^{\otimes k}(E_i)v_+=0$ for $i\le m$ as expected.
\end{proof}

Now from Proposition~\ref{pro:gyoja}, Theorem~\ref{thm:centra},
Corollary~\ref{cor:double}, the double centralizer theorem, and
Theorem~\ref{thm:final}, we have
\begin{thm}
Let $\lambda \vdash k$ be a partition in $H(m,n)$. Let $T$ be a
standard tableau of shape $\lambda$. Then $\uq$-submodule
$y_T(q)(V^{\otimes k})$ is isomorphic to the  irreducible
$\uu$-module $V(\lambda)$. Moreover as an $\hec\times \uu$
bimodule
$$
V^{\otimes k}\simeq \bigoplus_{\substack{\lambda \vdash k \\ \lambda
\in H(m,n)}} H^\lambda \otimes V(\lambda),
$$
where $H^\lambda$ is the irreducible $\hec$-module labeled by
$\lambda$.
\end{thm}

\bibliographystyle{plain}
\bibliography{super}

\begin{thebibliography}{10}

\bibitem{bkk}
G.~Benkart, S.~Kang, and M.~Kashiwara.
\newblock Crystal basis for the quantum superalgebra $\mathcal{U}_q(gl(m,n))$.
\newblock {\em Journal of Amer. Math. Soc.}, 13:295--331, 2000.

\bibitem{nova}
G.~Benkart and C.~Lee.
\newblock Stability in modules for general linear {L}ie superalgebras.
\newblock {\em Nova Journal of Algebra and Geometry}, 2(4):383--409, 1993.

\bibitem{br}
A.~Berele and A.~Regev.
\newblock Hook {Y}oung diagrams with applications to combinatorics, and to
  representations of {L}ie superalgebras.
\newblock {\em Adv. in Math.}, 64:118--175, 1987.

\bibitem{flv}
R.~Floreanini, D.~Leites, and L.~Vinet.
\newblock On the defining relations of quantum superalgebras.
\newblock {\em Lett. in Math. Physics}, 23:127--131, 1991.

\bibitem{frob}
F.~Frobenius.
\newblock {\"{U}}ber die {C}haraktere der symmetricschen {G}ruppe.
\newblock {\em Preuss. Akad. Wiss. Sitz.}, pages 516--534, 1900.
\newblock reprinted in Gesamelte Abhandlungen vol. 3, 148--166.

\bibitem{gyoja}
A.~Gyoja.
\newblock A $q$-analogue of {Y}oung symmetrizer.
\newblock {\em Osaka J. Math,}, 23:841--852, 1986.

\bibitem{jimbo}
M.~Jimbo.
\newblock A $q$-analog of $\mathcal{U}(gl(n+1))$, {H}ecke algebra, and the
  {Y}ang-{B}axter equation.
\newblock {\em Lett. Math. Phys.}, 11:247--252, 1986.

\bibitem{kt}
S.M. Khoroshkin and V.N. Tolstoy.
\newblock Universal ${R}$-matrix for quantized (super)algebras.
\newblock {\em Commun. Math. Phys.}, 141:599--617, 1991.

\bibitem{lr}
R.~Leduc and A.~Ram.
\newblock A {R}ibbon {H}opf algebra approach to the irreducible representations
  of centralizer algebras:{T}he {B}rauer {B}irman-{W}enzl, and type {A}
  {I}wahori-{H}ecke algebras.
\newblock {\em Adv. in Math.}, 125:1--94, 1997.

\bibitem{lusz81}
G.~Lusztig.
\newblock On a theorem of {B}enson and {C}urtis.
\newblock {\em J. Algebra}, 71:490--498, 1981.

\bibitem{scheu92}
M.~Scheunert.
\newblock Serre-type relations for special linear {L}ie superalgebras.
\newblock {\em Lett. Math. Phys.}, 24:173--181, 1992.

\bibitem{scheu93}
M.~Scheunert.
\newblock Presentation and $q$ deformation of {L}ie superalgebras.
\newblock {\em J. Math. Phys.}, 34:3780--3808, 1993.

\bibitem{schur1}
I.~Schur.
\newblock {\em {\"{U}}ber eine Klasse von Matrizen, die sich einer gegeben
  Matrix zuordenen lassen}.
\newblock PhD thesis, 1901.
\newblock reprinted in Gesamelte Abhandlungen vol. 1, 1--70.

\bibitem{schur2}
I.~Schur.
\newblock {\"{U}}ber die rationalen {D}arstellungen der allgemeinen linearen
  {G}ruppe.
\newblock {\em Preuss. Akad. Wiss. Sitz.}, pages 58--75, 1927.
\newblock reprinted in Gesamelte Abhandlungen vol. 3, 68--85.

\bibitem{young}
A.~Young.
\newblock On {Q}uantitative substitutional analysis i-ix, 1901--1952.
\newblock reprinted in {\it The Collected Papaers of Alfred Young 1873-1940},
  Mathematical Exposition No.{\bf 21} U. of Toronto Press,(1977).

\end{thebibliography}

\end{document}